\documentclass[preprint,12pt]{elsarticle}

\usepackage{amsmath,amsfonts,amssymb,bm,mathtools}
\usepackage{xcolor}
\usepackage{lineno}
\usepackage{linenofix}
\usepackage{subcaption}
\usepackage{booktabs}
\usepackage{siunitx}
\usepackage{algorithm,algpseudocode}
\usepackage{hyperref}
\usepackage[capitalize,nameinlink]{cleveref}

\hypersetup{colorlinks=true}

\begin{document}

\begin{frontmatter}

  \title{Physics-Informed Machine Learning with Conditional Karhunen-Lo\`eve Expansions}

  \author[PNNL]{A.M.~Tartakovsky\corref{cora}}%
  \ead{Alexandre.Tartakovsky@pnnl.gov}%

  \author[PNNL]{D.A.~Barajas-Solano}%
  \ead{David.Barajas-Solano@pnnl.gov}%

  \author[PNNL]{Q. He}%
  \ead{Qizhi.He@pnnl.gov}%

  \cortext[cora]{Corresponding author}%
  \address[PNNL]{Pacific Northwest National Laboratory, Richland, WA 99354}

  \begin{abstract}
    We present a new physics-informed machine learning approach for the inversion of PDE models with heterogeneous parameters. 
    In our approach, the space-dependent partially-observed parameters and states are approximated via Karhunen-Lo\`eve expansions (KLEs). 
    Each of these KLEs is then conditioned on their corresponding measurements, resulting in low-dimensional models of the parameters and states that resolve observed data. 
    Finally, the coefficients of the KLEs are estimated by minimizing the norm of the residual of the PDE model evaluated at a finite set of points in the computational domain, ensuring that the reconstructed parameters and states are consistent with both the observations and the PDE model to an arbitrary level of accuracy.
  
    In our approach, KLEs are constructed using the eigendecomposition of covariance models of spatial variability. 
    For the model parameters, we employ a parameterized covariance model calibrated on parameter observations; for the model states, the covariance is estimated from a number of forward simulations of the PDE model corresponding to realizations of the parameters drawn from their KLE. 
    We apply the proposed approach to identifying heterogeneous log-diffusion coefficients in diffusion equations from spatially sparse measurements of the log-diffusion coefficient and the solution of the diffusion equation.
    We find that the proposed approach compares favorably against state-of-the-art point estimates such as maximum a posteriori estimation and physics-informed neural networks.
  \end{abstract}

    
\end{frontmatter}

\newlength{\threesubht}
\newsavebox{\threesubbox}



\section{Introduction}
\label{sec:problem-formulation}

Parameter estimation is a critical step in modeling natural and engineered systems~\cite{vrugt2008inverse}. Here, we propose a new physics-informed machine learning method for estimating both parameters and states in systems described by differential equations. 
We consider the behavior of stationary physical systems modeled by PDEs over the simulation domain $D \subset \mathbb{R}^d$, $d \in [1, 3]$.
For simplicity, we assume that the system can be described by a single spatially heterogeneous scalar parameter $y \colon D \to \mathbb{R}$, one state variable $u \colon D \to \mathbb{R}$, and the stationary PDE problem $\mathcal{L}(u, y) = 0$, where $\mathcal{L}(\cdot, \cdot)$ denotes the governing equation and boundary conditions.
In this context, the ``forward'' problem is the problem of computing $u$ given $y$, and the ``inverse'' problem is the problem of estimating both $y$ and $u$ given measurements of $y$ and $u$.
In this work, we focus on the inverse problem with spatially sparse measurements of $y$ and $u$.

We assume that $N^u_{\mathrm{s}}$ measurements of $u$, $\{u_i\}^{N^u_{\mathrm{s}}}_{i = 1}$, are collected at spatial locations $\{ x^u_i \}^{N^u_{\mathrm{s}}}_{ i = 1}$.
Similarly, $N^y_{\mathrm{s}}$ measurements of $y$, $\{y_i \}^{N^y_{\mathrm{s}}}$ are collected at locations $\{ x^y_i \}^{N^y_{\mathrm{s}}}_{i = 1}$.
The observations are organized into the vector of observations $\mathbf{u}_{\mathrm{s}} = (u_1, \dots, u_{N^u_{\mathrm{s}}} )^{\top}$ and $\mathbf{y}_{\mathrm{s}} = (y_1, \dots, y_{N^y_{\mathrm{s}}} )^{\top}$, while the observation locations are organized into the observation matrices $X^u_{\mathrm{s}} = (x^u_1, \dots, x^u_{N^u_{\mathrm{s}}})$ and $X^y_{\mathrm{s}} = (x^y_1, \dots, x^y_{N^y_{\mathrm{s}}})$.
Finally, we assume that the observations are contaminated by normally distributed observation error, and we denote by $\Sigma_u$ and $\Sigma_y$ the error covariance matrices of the $u$ and $y$ observations, respectively.

The inverse problem can be defined as finding the functions $u$ and $y$ that minimize the discrepancy with respect to the observed data while satisfying the governing equations and boundary conditions \citep{golmohammadi2018exploiting,elsheikh2014efficient}, that is,
\begin{equation}
  \label{eq:model-inversion}
  \begin{aligned}
    \min_{u, y} \quad & \| u(X^u_{\mathrm{s}}) - \mathbf{u}_{\mathrm{s}} \|^2_{\Sigma_u} + \| y(X^y_{\mathrm{s}}) - \mathbf{y}_{\mathrm{s}} \|^2_{\Sigma_y}, \\
    \textrm{s.t.} \quad & \mathcal{L}(u, y) = 0.
  \end{aligned}
\end{equation}
where $\| \mathbf{v} \|_{\Sigma} \coloneqq \mathbf{v}^{\top} \Sigma^{-1} \mathbf{v}$ denotes the $\ell_2$ norm of the vector $\mathbf{v}$ weighted by the inverse of the covariance matrix $\Sigma$.

The problem of \cref{eq:model-inversion} is often solved numerically by discretizing the fields $y$ and $u$ and replacing the PDE constraint with its weak form corresponding to the discretization scheme.
Let $N$ denote the number of degrees of freedom of the discretization of the PDE problem.
In general, $N \gg N_s^u+N_s^y$, and the optimization problem of \cref{eq:model-inversion}) is ill-posed and requires regularization to have a unique solution \cite{engl1989convergence}.
The regularized problem reads
\begin{equation}
  \label{eq:map}
  \begin{aligned}
    \min_{u, y} \quad & \| u(X^u_{\mathrm{s}}) - \mathbf{u}_{\mathrm{s}} \|^2_{\Sigma_u} + \| y(X^y_{\mathrm{s}}) - \mathbf{y}_{\mathrm{s}} \|^2_{\Sigma_y} + \gamma \mathcal{R}(y), \\
    \textrm{s.t.} \quad & \mathcal{L}(u, y) = 0,
  \end{aligned}
\end{equation}
where $\mathcal{R}(\cdot)$ is a regularization penalty encoding regularization assumptions.
The regularization parameter $\gamma > 0$ controls the degree to which the discrepancy terms are minimized versus how much the regularization term is minimized.
In the context of Bayesian inference \cite{stuart-2010-inverse,ma2009efficient,marzouk2009dimensionality,burger-2014-maximum}, up to additive constants, the discrepancy terms are equivalent to the negative log likelihood of the observations, and $\gamma \mathcal{R}(y)$ is equivalent to the negative log prior density of $y$.
Therefore, the solution for $y$ of \cref{eq:map} is equivalent to the so-called \emph{maximum a posteriori} (MAP) estimate, a Bayesian point estimate defined as the largest mode of the posterior density of $y$ conditional on the observations.
Common choices for $\mathcal{R}(\cdot)$ include the so-called $H_1$ norm, $\| \nabla (\cdot) \|^2_2$, and total variation denoising (TVD), $\| \nabla (\cdot) \|_1$~\cite{barajassolano-2014-linear}.

Another approach to regularize the optimization problem of \cref{eq:model-inversion} is the pilot point method 
\cite{doherty2003ground,alcolea2006pilot,ramarao1995pilot}. This method consist of parametrizing $y$ in terms of its value at a set of so-called ``pilot points''.
Everywhere else in $D$, $y$ is evaluating by regressing $y$ measurements and the pilot point values using, e.g., Gaussian Process regression (also known as ``kriging'')~\cite{stein-1999-interpolation,williams2006gaussian,cressie1990origins,gunes2006gappy,cressie-2015-geostatistics}.
The value of $y$ at the pilot point locations are estimated from the minimization problem of \cref{eq:model-inversion}).


Bayesian methods, such as Ensemble Kalman Filter (EnKF)~\citep{chen2013application,evensen2009data,camporese2009ensemble,stuart-2017-analysis,xu-simultaneous-2018} and cokriging~\cite{mclaughlin1995recent,yang-2019-physics}, are commonly used for approximately solving the inversion problem~\eqref{eq:model-inversion}.
Following stochastic approach to modeling flow and transport \cite{dagan2005subsurface}, EnKF and cokriging treat $y$ and $u$ as the random fields $y(x,\omega) = \mathbb{E} \left [ y(x,\omega) \right ] + y'(x,\omega)$ and $u(x,\omega) = \mathbb{E} \left [ u(x,\omega) \right ] + u'(x,\omega)$ with expectations $\bar{y}(x) \coloneqq \mathbb{E} \left [ y(x,\omega) \right ] = $ and $\bar{u}(x) \coloneqq \mathbb{E} \left [ u(x,\omega) \right ]$, and zero-mean fluctuations $u'(x,\omega)$ and $y'(x,\omega)$.
The parameter estimate is computed using a cokriging update rule of the form
\begin{multline}
  \label{eq:EnKF}
  y^{\mathrm{EnKF}}(x) = \bar{y}(x) + C_{yu} \left ( x, X^u_{\mathrm{s}} \right) \left [ C_u(X^u_{\mathrm{s}}, X^u_{\mathrm{s}}) \right ]^{-1} \left [ \mathbf{u}_{\mathrm{s}} - \bar{u}(X^u_{\mathrm{s}}) \right ]\\
  + C_y \left ( x, X^y_{\mathrm{s}} \right) \left [ C_y(X^y_{\mathrm{s}}, X^y_{\mathrm{s}}) \right ]^{-1} \left [ \mathbf{y}_{\mathrm{s}} - \bar{y}(X^y_{\mathrm{s}}) \right ],
\end{multline}
where $C_y$ and $C_u$ denote the covariances of the $y$ and $u$ random fields, respectively, and $C_{yu}$ denotes the $y$-$u$ cross-covariance.
These covariances are evaluated in practice using sample-based estimates.
Inversion schemes of the form of~\cref{eq:EnKF} are straightforward to implement and do not require directly solving a minimization problem.
Nevertheless, the resulting estimate $y^{\mathrm{EnKF}}$ is not consistent with both data and physics; that is, the solution $u$ of $\mathcal{L}(u, y^{\mathrm{EnKF}}) = 0$ does not match the observed at $\mathbf{u}_s$.
Fully Bayesian methods~\cite{stuart-2010-inverse} address this inconsistency but often incur in significant computational effort, although significant advances have been made in recent years to address computational cost~\cite{barajassolano-2019-approximate,beskos-2017-geometric}.

Machine Learning (ML) methods have arisen in recent years as popular approaches for scientific applications.
In general, ML methods require a large amount of data and therefore are not feasible for parameter estimation with sparse measurements.
To address this challenge, a physics-informed neural networks (PINNs)
\cite{raissi2017physicsPart2,raissi2017physicsPart1,raissi2018deep} was extended for solving the inverse problem of \cref{eq:model-inversion}~\cite{tartakovsky2018learning}.
In this method, both $y$ and $u$ are represented with feed-forward deep neural networks as $u(x) \approx \hat{u}(x; \bm{\theta})$ and $y(x) \approx \hat{y}(x; \bm{\gamma})$, where $\bm{\theta}$ and $\bm{\gamma}$ denote the vectors of neural network weights.
Next, a ``residual'' network is defined as
\begin{equation}
  \hat{r} (x; \bm{\theta}, \bm{\gamma}) =  \mathcal{L} \left (\hat{u}(x; \bm{\theta}), \hat{y}(x; \bm{\gamma}) \right ),
\end{equation}
where differentiation with respect to $x$ is performed using automatic differentiation. 
These three networks are trained jointly by minimizing the loss function
\begin{equation}
  \label{eq:loss_fn_linear} 
  \min_{\bm{\theta}, \bm{\gamma}}  
  \left \| \hat{u}(X^u_{\mathrm{s}}; \bm{\theta}) - \mathbf{u}_{\mathrm{s}} \right \|^2_{\Sigma_u} + \left \| \hat{y}(X^y_{\mathrm{s}}; \bm{\gamma}) - \mathbf{y}_{\mathrm{s}} \right \|^2_{\Sigma_y} 
  + \rho \left \| \hat{r} (X^{r}; \bm{\theta}, \bm{\gamma}) \right \|^2_2,
\end{equation}
where the residual network is evaluated at certain ``residual'' points $\{ x^r_i \in D \}^{N_r}_{i = 1}$, organized into the matrix $X^r = (x^r_1, \dots, x^r_{N_r})$.
In this approach, the PDE constraint in~\cref{eq:model-inversion} is replaced with a weaker constraint on the residuals of $\mathcal{L}(u, y)$; therefore, the estimated fields only approximately satisfy the physics.
The advantage of PINNs is that it does not require discretizing the governing PDE for solving inverse problems.

Here, we propose a new physics-informed ML method for inverse problems based on conditional Karhunen Lo\`{e}ve expansions (cKLEs)~\cite{tipireddy2019conditional}.
In our approach, we model the fields $y$ and $u$ as realizations of Gaussian random fields $\hat{y}^c$ and $\hat{u}^c$ conditioned on observed data.
These Gaussian random fields encode the spatial correlation structure of the fields $y$ and $u$; for the $u$ variable, the corresponding random field satisfies both the data and governing PDE problem $\mathcal{L}(u, y) = 0$.
For both random fields, we compute their cKLEs, which allow us to parametrize their realizations in terms of so-called KL coefficients.
These KL coefficients are then estimated by solving a regularized form of~\cref{eq:model-inversion}.
We refer to cKLEs trained in this manner as ``physics-informed cKLEs'', or PICKLEs.

Similarly to PINNs, PICKLEs are trained to satisfy the governing equation $\mathcal{L}(u, y) = 0$ by penalizing the norm of a vector of residuals.
Unlike deep neural networks, the KLE of a field enforces its spatial correlation structure and acts as a regularizer.
Our results indicate that if the correlation structure of the underlying fields to be estimated is known or can be well estimated from observation data, then the PICKLE method for inverse problems leads to more accurate parameter estimates than such state-of-the-art inversion approaches as MAP estimation, or PINNs.

The remainder of this manuscript is structured as follows.
In~\cref{sec:ckle}, we introduce cKLEs.
We describe our algorithm for inverse problems based on PICKLEs in~\cref{sec:ckli}.
Finally, in~\cref{sec:numerical}, we apply the PICKLE method for the inverse problem of estimating the heterogeneous log-diffusion coefficient of the diffusion equation from sparse measurements of the log-diffusion coefficient and the solution of the diffusion equation.
PICKLE estimates are found to compare favorably against MAP and PINNs estimates.

\section{Conditional Karhunen Lo\`{e}ve expansions}
\label{sec:ckle}

Karhunen Lo\`{e}ve expansions (KLEs)~\cite{huang2001convergence} are used for representing random fields in terms of linear combinations of uncorrelated random variables.
In this work, we employ KLEs as parameterized, deterministic, representations of $u$ and $y$.
Specifically, we treat partially known $u$ and $y$ as realizations random fields $\hat{u}^c : D \times \Omega \to \mathbb{R}$ and $\hat{y}^c : D \times \Omega \to \mathbb{R}$ (where $\Omega$ is the corresponding random outcome space) conditioned on observed data.
Next, we compute the KLEs of these fields, which we use to parametrize their realizations.
We refer to these KLEs as \emph{conditional} KLEs, or cKLEs, as by construction they resolve observed data, i.e., at the observation locations the cKLE mean is equal to the field's observation and the cKLE variance is equal to the observation error variance.

In this section we discuss the construction of the cKLEs.
The selection of the Gaussian random field models $\hat{u}^c$ and $\hat{y}^c$ is discussed in \cref{sec:ckli-ckle}.

To introduce cKLEs, we consider a Gaussian random field $z \colon D \times \Omega \to \mathbb{R}$ with the expectation and covariance function, respectively,%
\begin{equation*}
  \bar{z}(x) \coloneqq \mathbb{E} \left [ z (x,\omega) \right ], \quad C(x, x') \coloneqq \mathbb{E} \left \{ [ z(x,\omega) - \bar{z}(x) ] [z(x',\omega) - \bar{z}(x') ] \right \}.
\end{equation*}
Next, we assume that a number of noisy spatial observations of $z$ are available, and similar to \cref{sec:problem-formulation}, these observations and the observation locations are organized into the vector $\mathbf{z}_{\mathrm{s}}$ and the matrix $X_{\mathrm{s}}$, respectively.
Furthermore, we denote by $C_s \coloneqq C(X_{\mathrm{s}}, X_{\mathrm{s}}) + \Sigma$ the covariance matrix of the observations, where $\Sigma$ is the covariance matrix of observation errors.
Employing Gaussian process regression (GPR)~\cite{williams2006gaussian}, we find that the conditional Gaussian process (GP) $z^c(x,\omega) \coloneqq z(x, \omega) \mid \left (\mathbf{z}_{\mathrm{s}}, X_{\mathrm{s}} \right )$ has the conditional mean and covariance kernels
\begin{align}
  \label{eq:gpr-mean}
  \bar{z}^c(x) &= \bar{z}(x) + C(x, X) C^{-1}_{\mathrm{s}} \left [ \mathbf{z}_{\mathrm{s}} - \bar{z}(X_{\mathrm{s}}) \right ],\\
  \label{eq:gpr-var}
  C^c(x, x') &= C(x, x') - C(x, X) C^{-1}_{\mathrm{s}} C(X, x'),
\end{align}
where the superindex $c$ stands for ``conditional'' on observations.

The cKLE of $z$, $z(x, \bm{\xi}(\omega)) = z(x, \omega)$, reads
\begin{equation}
  \label{eq:ckle}
  z^c \left ( x, \bm{\xi}(\omega) \right ) = \bar{z}^c(x) + \sum^{\infty}_{i = 1} \phi_i(x) \sqrt{\lambda_i} \xi_i(\omega) ,
\end{equation}
%
%
where $\bm{\xi}(\omega) = (\xi_1(\omega), \xi_2(\omega), \cdots )^{\top}$ is a vector of zero-mean, independent, identically-distributed standard Gaussian random variables, and the eigenpairs $\{ \phi_i(x), \lambda_i \}^{\infty}_{i = 1}$ are the solutions to the eigenvalue problem
\begin{equation*}
  \int_D C^c(x, x') \phi(x') \, \mathrm{d} x' = \lambda \phi(x).
\end{equation*}
%
The sequence of eigenfunctions forms an orthonormal basis on $L_2(D)$.

As the sum in \cref{eq:ckle} is infinite, the cKLE in this form is not directly amenable to numerical calculations.
Instead, in this work we will truncate cKLEs to a finite number of terms.
For random fields with non-trivial correlation structures (i.e., the non-zero correlation length), the eigenspectrum (i.e., the sequence of eigenvalues $\lambda_i$) decays towards zero for increasing $i$.
This, together with the Mercer theorem, justifies the truncation of the KLE to a finite number of terms~\cite{spanos1989stochastic}.
By the Mercer theorem, the KLE truncated to $M$ terms,
\begin{equation}
  \label{eq:cKLE-truncated}
  z^c_M \left (x, \bm{\xi}_M(\omega) \right ) = \bar{z}^c(x) + \sum^M_{i = 1} \phi_i(x) \sqrt{\lambda_i} \xi_i (\omega)
\end{equation}
converges to $z^c(x, \omega)$ in the $L_2$ sense for increasing $M$, that is,
\begin{equation*}
  \mathbb{E} \left \{ \left [ z^c(x,\omega) - z^c_M \left ( x, \bm{\xi}(\omega) \right ) \right ]^2 \right \} = \sum^{\infty}_{i = M + 1} \lambda_i \phi^2_i(x) .
\end{equation*}
This statement of convergence provides a means for selecting $M$ a priori in the context of uncertainty quantification.
By the orthonormality of the basis, it follows that the bulk variance and the mean-square truncation error are given by 
\begin{equation*}
  \int_D \operatorname{Var} z^c(x) \, \mathrm{d} x = \sum^{\infty}_{i = 1} \lambda_i
\end{equation*}
and  
\begin{equation}
  \label{eq:mse-truncation}
  \int_D \mathbb{E} \left \{ \left [ z^c(x, \omega) - z^c_M(x, \bm{\xi}(\omega)) \right ]^2 \right \} \, \mathrm{d} x = \sum^{\infty}_{i = M + 1} \lambda_i,
\end{equation}
respectively.
Therefore, $M$ is commonly chosen based on either of the following relative and absolute conditions
\begin{equation}
  \label{eq:truncation}
  \sum^{\infty}_{i = M + 1} \lambda_i \leq \text{rtol} \int_D \operatorname{Var} z(x) \, \mathrm{d} x, \quad \sum^{\infty}_{i = M + 1} \lambda_i \leq \text{atol},
\end{equation}
for certain relative and absolute tolerances $\text{rtol}$ and $\text{atol}$, respectively.
We must note that~\cref{eq:mse-truncation} is a statement about the bulk squared truncation error averaged over all realizations of $z^c$, and does not provide a bound for the bulk squared truncation error for any given realization.

For the sake of brevity, we organize the sequences of eigenvalues, eigenvectors, and random variables into the vector of functions
\begin{equation}
  \label{eq:dot-prod-components}
  \bm{\psi}(x) = \left [ \sqrt{\lambda_1} \psi_1(x), \cdots, \sqrt{\lambda_M} \psi_M(x) \right ]^{\top}, \quad \bm{\xi}(\omega) = (\xi_1(\omega), \cdots, \xi_M(\omega) )^{\top}
\end{equation}
so that the truncated cKLE~\cref{eq:cKLE-truncated} can be rewritten in dot product form as
\begin{equation}
  \label{eq:truncated-KLE-dot-prod}
  z^c_M \left (x, \bm{\xi}(\omega) \right ) = \bar{z}(x) + \bm{\psi}^{\top}(x) \bm{\xi}(\omega).
\end{equation}

If we treat the $\xi_i$s in the sequence $\{\xi_i\}^{M}_{i = 1}$ not as random variables but as expansion coefficients, we can understand the cKLE as a parameterized representation of functions that satisfy up to measurement error the observed data $(\mathbf{z}_{\mathrm{s}}, X_{\mathrm{s}})$.
In this context, we refer to the $\xi_i$s as the cKLE ``coefficients''.
Estimating a certain function that satisfies the observed data is then a matter of estimating the cKLE coefficients.
We will employ this interpretation of cKLEs to construct our parameter estimation approach in the following section.

\section{cKLE-based inversion}
\label{sec:ckli}

In this section, we present the PICKLE method for parameter estimation.
In~\cref{sec:ckli-ckle}, we describe the selection of the Gaussian random fields used to construct the cKLEs of $y$ and $u$.
In~\cref{sec:ckli-minimization} we describe how we train these cKLEs subject to a PDE constraint.

\subsection{Constructing cKLEs of $y$ and $u$}
\label{sec:ckli-ckle}

To construct the cKLEs of $y$ and $u$, we first construct conditional GPs $\hat{y}^c \colon D \times \Omega \to \mathbb{R}$ and $\hat{u}^c \colon D \to \Omega \to \mathbb{R}$.
Specifically, we select the (unconditional) mean and covariance kernel of these GPs so that they encode the spatial correlation structure of the fields to be estimated.
Once the unconditional mean and covariance kernel are selected, the conditional random fields $\hat{y}^c$ and $\hat{u}^c$ are obtained by conditioning on observation data using~\cref{eq:gpr-mean,eq:gpr-var}.


\subsubsection{cKLE of $y$}
\label{sec:ckli-ckle-y}

For $\hat{y}^c$, we set the unconditional mean to zero and select the unconditional covariance kernel from a parameterized family of covariance kernels $C^y(\cdot, \cdot \mid \bm{\theta})$ such as the Mat\'{e}rn, exponential, or square exponential (i.e., Gaussian) kernels.
The parameters $\bm{\theta}$ of the kernel are estimated from the observation data $(\mathbf{y}_{\mathrm{s}}, X^y_{\mathrm{s}})$ via marginal likelihood maximization or leave-one-out cross-validation as is commonly done in GPR.
We justify this GPR-based approach by noting that it is commonly used in geophysics, under the name of kriging, for estimating spatially heterogeneous geophysical parameters from sparse observations.

Once the unconditional covariance kernel is selected and the conditional mean and covariance of $\hat{y}^c$ are evaluated (via~\cref{eq:gpr-mean,eq:gpr-var}), we construct the cKLE  of $\hat{y}^c$ truncated to $N_{\xi}$ terms of the form of \cref{eq:truncated-KLE-dot-prod}, namely,
\begin{equation}
  \label{eq:cKLE-y-random}
  \hat{y}^c (x, \bm{\xi}) = \bar{y}^c(x) + \bm{\psi}^{\top}_y (x) \bm{\xi}, \quad \bm{\xi} \sim \mathcal{N}(\bm{0}, \mathbf{I}_{N_{\xi}}),
\end{equation}
where $\mathcal{N}(0, \mathbf{I}_{N_{\xi}})$ is the multivariate normal distribution and $\mathbf{I}_{N_{\xi}}$ is the $N_{\xi} \times N_{\xi}$ identity matrix.

\subsubsection{cKLE of $u$}
\label{sec:ckli-ckle-u}

For $\hat{u}^c$, the data-driven GPR-based strategy is inadequate as samples from common parameterized Gaussian process models are not guaranteed to satisfy the governing equations and boundary conditions.
Therefore, in this work we employ a Monte Carlo simulation-based method for computing the unconditional mean and covariance for $u$.
We construct an ensemble of $N_{\mathrm{ens}}$ realizations of $\hat{y}^c$, $\{ y^{(i)} \}^{N_{\mathrm{ens}}}_{i = 1}$ by sampling $\bm{\xi}^{(i)}$ from $\mathcal{N}(0, \mathbf{I}_{N_{\xi}})$ and, then, evaluating the cKLE model of $\hat{y}^c$, \cref{eq:cKLE-y-random}, with $\bm{\xi} = \bm{\xi}^{(i)}$, that is,
\begin{equation}
  \label{eq:uc-ysample}
  y^{(i)}(x) = \bar{y}^c(x) + \bm{\psi}^{\top}_y (x) \bm{\xi}^{(i)}, \quad \bm{\xi}^{(i)} \sim \mathcal{N}(\bm{0}, \mathbf{I}_{N_{\xi}}).
\end{equation}
For each member of the ensemble $\{ y^{(i)} \}^{N_{\mathrm{ens}}}_{i = 1}$, we calculate $u^{(i)}$ by solving the PDE problem $\mathcal{L}(u^{(i)}, y^{(i)}) = 0$, thus obtaining the ensemble of $u$ fields, $\{ u^{(i)} \}^{N_{\mathrm{ens}}}_{i = 1}$.
The unconditional mean and covariance for $\hat{u}^c$, $\overline{u}$ and $C_u(x,x')$, are then computed as the ensemble estimates
\begin{align}
  \label{eq:ens_mean}
  \overline{u}(x) &= \frac{1}{N_{\mathrm{ens}}} \sum_{i=1}^{N_{\mathrm{ens}}} u^{(i)}(x), \\
  \label{eq:ens_covariance}
  C_u(x,x') &= \frac{1}{N_{\mathrm{ens}} - 1} \sum_{i=1}^{N_{\mathrm{ens}}} \left [ u^{(i)}(x) - \overline{u}(x) \right ] \left [ u^{(i)}(x') - \overline{u}(x') \right ].
\end{align}
This procedure is summarized in~\cref{alg:uc}.

\begin{algorithm}
  \caption{Sampling-based covariance model for $u$}
  \label{alg:uc}
  \begin{algorithmic}[1]
    \Require $X^u_{\mathrm{s}}$, $\mathbf{u}_{\mathrm{s}}$, $N_{\mathrm{ens}}$
    \For{$i \gets 1, N_{\mathrm{ens}}$}
    \State Generate $y^{(i)}$ via~\eqref{eq:uc-ysample}
    \State Compute $u^{(i)}$ by solving $\mathcal{L}(u^{(i)}, y^{(i)}) = 0$
    \EndFor
    \State Compute ensemble mean and covariance of $\{u^{(i)}\}$, $\overline{u}$ and $C_u(x,x')$, using \cref{eq:ens_mean,eq:ens_covariance}
    \State Compute conditional mean and covariance of $\hat{u}^c$, $\overline{u}^c$ and $C_u^c(x,x')$, using \cref{eq:gpr-mean,eq:gpr-var}
  \end{algorithmic}
\end{algorithm}

Once the unconditional mean and covariance kernels are estimated, the conditional mean and covariance of $\hat{u}^c$ are calculated using~\cref{eq:gpr-mean,eq:gpr-var} and the cKLE model for $\hat{u}^c$ is constructed in the form of~\cref{eq:truncated-KLE-dot-prod}, namely,
\begin{equation}
  \label{eq:ckle-u-random}
  \hat{u}^c(x, \bm{\eta}) = \bar{u}^c(x) + \bm{\psi}^{\top}_u (x) \bm{\eta}.
\end{equation}


The ensemble covariance estimate requires some discussion.
We require $N_{\mathrm{ens}} > N^u_{\mathrm{s}}$ so that the rank of the unconditional covariance is larger than $N^u_{\mathrm{s}}$ and the conditional covariance is not trivial.
This limitation can be avoided by instead employing shrinkage estimators, which are known to be robust for small number of ensemble elements~\cite{chen_2009_shrinkage}.
In~\cite{yang-2019-physics}, the computational cost of the Monte Carlo simulations for estimating the unconditional covariance of $u$ was reduced by using the Multilevel Monte Carlo method~\cite{giles-2015-multilevel}.
For small unconditional variance of $y$, the moment equation method can be used to derive a system of deterministic equations for the unconditional covariance of $u$~\cite{Jarman2013JUQ,TARTAKOVSKY2003182}.
Then, the unconditional covariance of $u$ can be found by solving these equations numerically.  

\subsection{The PICKLE method for inverse problems}
\label{sec:ckli-minimization}

In this section we describe the proposed PICKLE method for inverse problems.
Similarly to PINN, in PICKLE we replace the PDE constraint in~\cref{eq:model-inversion} with a penalty on the norm of the vector of residuals,
\begin{equation*}
  \mathbf{r} [ u, y ] \coloneqq \left [ \mathcal{L} (u, y) \mid_{x = x^r_1}, \dots, \mathcal{L} (u, y) \mid_{x = x^r_{N_r}} \right ]^{\top},
\end{equation*}
where each component of the vector corresponds to the residual of the PDE problem evaluated at the $i$th ``residual'' point of the sequence $\{ x^r_i \in D \}^{N_r}_{i = 1}$.

The constraint on the residuals is added as a penalty term into the objective function, leading to the minimization problem
\begin{equation}
  \label{eq:model-inversion-penalized}
  \min_{u, y} \quad \| u(X^u_{\mathrm{s}}) - \mathbf{u}_{\mathrm{s}} \|^2_{\Sigma_u} + \| y(X^y_{\mathrm{s}}) - \mathbf{y}_{\mathrm{s}} \|^2_{\Sigma_y} + \rho \| \mathbf{r} [ u, y ] \|^2_2,
\end{equation}
where $\rho > 0$ is a penalty parameter.

We now proceed to introduce the cKLE models for $y$ and $u$.
Namely, we interpret the cKLEs \cref{eq:cKLE-y-random} and \cref{eq:ckle-u-random} as representations of functions parameterized by the vectors of cKLE coefficients $\bm{\xi}$ and $\bm{\eta}$, leading to the deterministic cKLE models
\begin{align}
  \label{eq:ckle-y}
  y^c(x, \bm{\eta}) &= \bar{y}^c(x) + \bm{\psi}^{\top}_y (x) \bm{\xi},\\
  \label{eq:ckle-u}
  u^c(x, \bm{\eta}) &= \bar{u}^c(x) + \bm{\psi}^{\top}_u (x) \bm{\eta}.
\end{align}
Substituting \cref{eq:ckle-y,eq:ckle-u} into~\cref{eq:model-inversion-penalized}, we obtain the following minimization problem in terms of the cKLE parameters:
\begin{equation*}
  \min_{\bm{\xi}, \bm{\eta}} \quad \| u^c(X^u_{\mathrm{s}}, \bm{\eta}) - \mathbf{u}_{\mathrm{s}} \|^2_{\Sigma_u} + \| y^c(X^y_{\mathrm{s}}, \bm{\xi}) - \mathbf{y}_{\mathrm{s}} \|^2_{\Sigma_y} + \rho \left \| \mathbf{r} [ u^c(\cdot, \bm{\eta}), y^c(\cdot, \bm{\xi}) ] \right \|^2_2.
\end{equation*}
By construction, the cKLE models minimize the discrepancy terms.
This leaves only the penalty term, so that the coefficient $\rho$ can be dropped.

It remains to regularize the problem.
In this work we choose to penalize the $\ell_2$-norm of the vectors of cKLE parameters, resulting in the final PICKLE minimization problem
\begin{equation}
  \label{eq:ckli}
  \min_{\bm{\xi}, \bm{\eta}} \quad \left \| \mathbf{r} [ u^c(\cdot, \bm{\eta}), y^c(\cdot, \bm{\xi}) ] \right \|^2_2 + \gamma \left ( \left \| \bm{\xi} \right \|^2_2 + \left \| \bm{\eta} \right \|^2_2 \right ),
\end{equation}
where $\gamma > 0$ is a regularization penalty.
Substituting $\bm{\xi}$ and $\bm{\eta}$, estimated from \cref{eq:ckli}, into~\cref{eq:ckle-y} and \cref{eq:ckle-u} provides the PICKLE estimates of the $y$ and $u$ fields.
The proposed model inversion algorithm is summarized in~\cref{alg:ckli}.

\begin{algorithm}
  \caption{cKLE-based inversion}
  \label{alg:ckli}
  \begin{algorithmic}[1]
    \Require $X^y_{\mathrm{s}}$, $\mathbf{y}_{\mathrm{s}}$, $X^u_{\mathrm{s}}$, $\mathbf{u}_{\mathrm{s}}$, $C^y(\cdot, \cdot \mid \bm{\theta})$, $N_{\xi}$, $N_{\eta}$, $N_{\mathrm{ens}}$
    \State Estimate $\bm{\theta}$ via GPR model selection
    \State Compute conditional mean and covariance of $y^c$ using \cref{eq:gpr-mean,eq:gpr-var}
    \State Calculate KLE of $y^c$
    \State Calculate cKLE model, \cref{eq:ckle-y}, truncated to $N_{\xi}$ terms
    \State Compute conditional mean and covariance of $u^c$ using \cref{alg:uc}
    \State Calculate KLE of $u^c$
    \State Calculate cKLE model, \cref{eq:ckle-u}, truncated to $N_{\eta}$ terms
    \State Estimate $\bm{\xi}$ and $\bm{\eta}$ via~\cref{eq:ckli}
    \State Compute $y$ and $u$ from estimated $\bm{\xi}$ and $\bm{\eta}$ using~\cref{eq:ckle-y} and \cref{eq:ckle-u}
  \end{algorithmic}
\end{algorithm}

\subsection{Computational cost}
\label{sec:ckli-discussion}

Common iterative, gradient-based approaches to the solution of the PDE-constrained optimization problem of \cref{eq:model-inversion} aim to minimize the objective function with respect to $y$, with $u$ given explicitly at every iteration of the procedure as the solution of the PDE constraint, $\mathcal{L}(u, y) = 0$, for given $y$.
The gradient of the objective function with respect to $y$ is then found by the application of the chain rule and the adjoint method, e.g., see~\cite{zhang_2017_adjoint}.
Such approaches require solving the PDE constraint at every step of the iteration process.
In contrast, in PICKLE there is no need to solve the governing PDE.  
Instead, our approach requires only evaluating the norm of the vector of residuals and its gradient with respect to the cKLE coefficients.

The calculation of the residuals' norm gradient deserves special consideration.
One can consider a strong or weak form of the PDE residual. The strong form of the PDE residual requires evaluating the spatial derivatives of the cKLEs of $y$ and $u$, which in turn requires obtaining the cKLE in terms of closed form functions.
While the eigenproblem for the cKLE cannot be exactly solved in closed form in general, closed-form approximations in terms of orthogonal polynomials (e.g. Chebyshev polynomials) can be obtained (e.g., \cite{sraj_2016_coordinate,liu_2017_chebyshev}).
The benefit of having the closed-form cKLEs is that the norm of residuals and its gradients can be evaluated programatically using automatic differentiation of the composition of the residual and the cKLEs.
In this work, we consider the residual of a weak form of the PDE constraint.
In this case, it suffices to solve the eigenproblems and compute the cKLEs of $y$ and $u$ on the discretized grids corresponding to the weak form of the PDE problem.
In \cref{sec:numerical} we discuss the FV approximation of the PDE problem employed for the numerical experiments presented in this work.

The three factors that chiefly control the computational cost of the PICKLE approach are (i) the number of samples in the ensemble $\{u^{(i)}\}$, $N_{\mathrm{ens}}$, (ii) the number of cKLE parameters, (iii) and the size of the vector of residuals.
In the following, we discuss these sources of computational cost one-by-one.

(i) In PICKLE, the governing PDE is solved $N_{\mathrm{ens}}$ times, a number specified a priori. In comparison, traditional gradient-based approaches to the solution of the PDE-constrained optimization problem require a number of solutions of the PDE constraint that cannot be controlled a priori.
Therefore, the computational cost of solving complex physics problems cannot in general be controlled a priori for such approaches. Also, in PICKLE, each of $N_{\mathrm{ens}}$ realizations can be run independently. Therefore,  PICKLE is trivially parallelizable, which can dramatically reduce the computational time associated with this cost.  

(ii) As we show in~\cref{sec:numerical}, the cKLEs allow us to represent the $y$ and $u$ fields in terms of a relatively low number of KLE parameters, which makes it possible to to tackle high-dimensional problems; specifically, accurate solutions can be obtained with a number of KLE parameters significantly less than the number necessary to represent the unconditional $y$ and $u$ fields.
The accuracy PICKLE strongly depends on the expressive capacity of truncated cKLEs, which is known to be limited for fields with sharp gradients or discontinuities due to the Gibbs phenomenon.
Nevertheless, piecewise-continuous $y$ fields can be treated with our approach by introducing latent fields, as described in \cref{sec:numerical-piecewise}.

(iii) With regards to the vector of residuals, the proposed inversion approach provides significant flexibility for the choice of residuals.
For the numerical experiments shown in~\cref{sec:numerical}, we employ the residuals of the finite volume (FV) discretization of the PDE constraint evaluated at a subset of the FV elements of the discretization.
The residuals' vector size of can be adjusted to reduce the computational cost of the inverse problem solution.

\section{Numerical experiments}
\label{sec:numerical}

In this section, we use PICKLE for estimating the heterogeneous diffusion coefficient of the elliptic diffusion equation.
Specifically, we consider the PDE problem
\begin{equation}
  \label{eq:numerical-pde}
  \begin{aligned}
    \nabla \cdot \left [ e^{y(x)} \nabla u(x) \right ] &= 0, && x \in D \coloneqq [0, 1]^2,\\
    u(x) &= 1, && x_1 = 0,\\
    u(x) &= 0, && x_1 = 1,\\
    e^{y(x)} \frac{\partial u(x)}{\partial x_2} &= 0, && x_2 = \{0, 1\},
  \end{aligned}
\end{equation}
where $u(x)$ is the PDE solution and $y(x)$ is the log-diffusion coefficient.
Among other problems, this equation describes saturated flow in heterogeneous porous media \citep{bear2013dynamics}. 
Our goal is to estimate the spatial distribution of $y$ from noiseless sparse observations of $y$ and $u$.

In our numerical experiments, we discretize the simulation domain into a uniform grid of $32 \times 32$ rectangular elements, for a total of $1024$ elements.
The PDE problem~\cref{eq:numerical-pde} is then discretized employing a cell-centered FV scheme and the two-point flux approximation.

The PICKLE method is implemented in Scientific Python~\citep{oliphant_2007_computing}, and all numerical experiments are executed on a Intel Xeon W-2135 workstation employing \textsf{GNU Parallel}~\citep{tange_2011_parallel}. 
To evaluate the accuracy of the PICKLE method, we compare the reconstructed log-diffusion field to the reference field employed to generate the synthetic observations.
Furthermore, we compare the reconstructed field against the MAP estimate with $H_1$ regularization, a commonly used PDE-constrained optimization-based inversion approach, and the PINN method~\citep{tartakovsky2018learning}.
For both MAP and PICKLE, we use the regularization parameter $\gamma = \num{1e-6}$.
The PINN method does not employ regularization, other than the regularization provided by physics constraints.
The PINN implementation details are given in~\cref{sec:numerical-continuous}. 

\subsection{Continuous diffusion field}
\label{sec:numerical-continuous}

\begin{table}
  \centering
    \caption{Properties of synthetic reference log-diffusion fields, and PICKLE estimation parameters}
  \begin{tabular}{lccccccc}
    & $\lambda$ & $\sigma$ & $N_{\xi}$ & $N_{\eta}$ & $N^y_{\mathrm{s}}$ & $N^u_{\mathrm{s}}$ \\
    \midrule
    Gaussian & 0.2 & 1.0 & 100 & 100 & 50 & 50 \\
    Mat\'{e}rn $\nu = 5/2$ & 0.2 & 1.0 & 100 & 100 & 50 & 50 \\
    Mat\'{e}rn $\nu = 3/2$ & 0.1 & 1.0 & 300 & 200 & 200 & 50 \\
  \end{tabular}
  \label{tbl:numerical-continuous-params}
\end{table}

\begin{figure}[tbhp]
  \centering
  \begin{subfigure}{0.32\textwidth}
    \includegraphics[width=\textwidth]{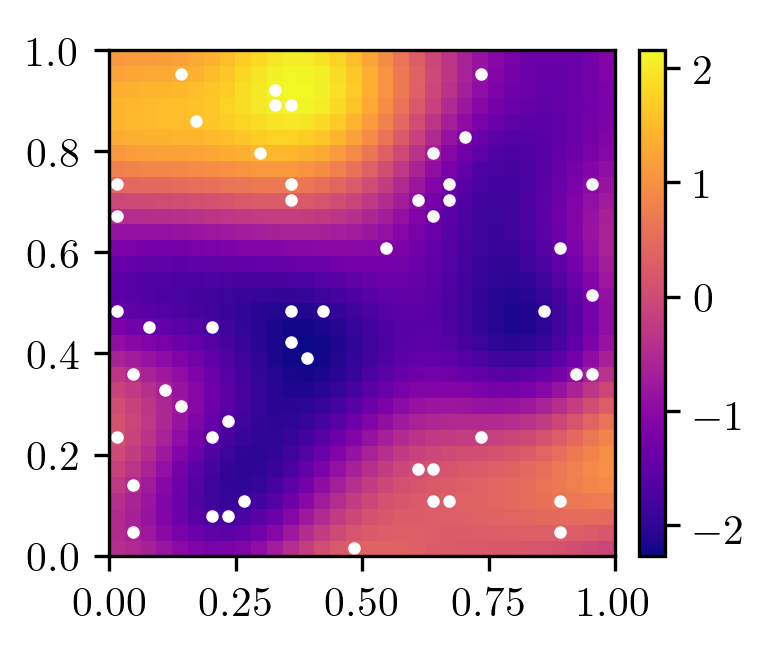}
    \caption{$y(x)$, Gaussian kernel}
  \end{subfigure}
  \begin{subfigure}{0.32\textwidth}
    \includegraphics[width=\textwidth]{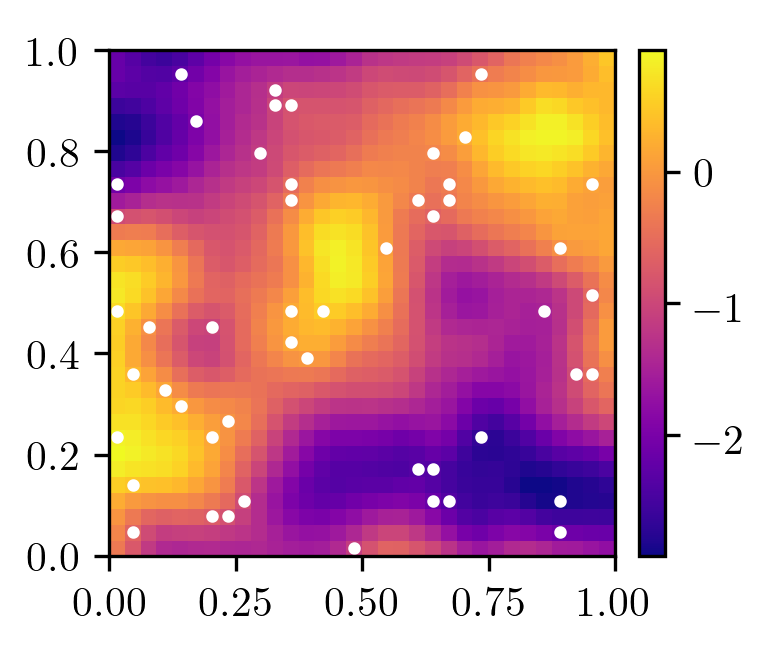}
    \caption{$y(x)$, Mat\'{e}rn $\nu = 5 / 2$}
  \end{subfigure}
  \begin{subfigure}{0.32\textwidth}
    \includegraphics[width=\textwidth]{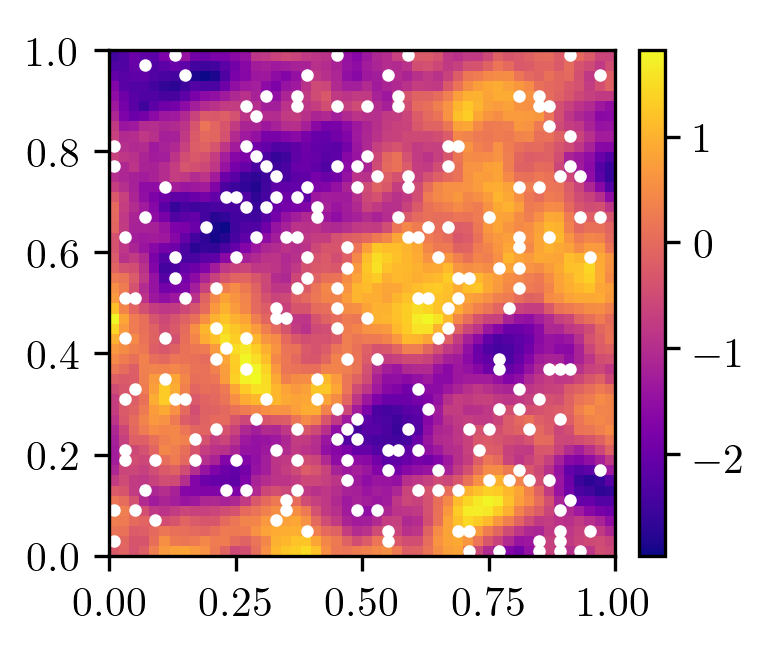}
    \caption{$y(x)$, Mat\'{e}rn $\nu = 3 / 2$}
  \end{subfigure}\\
  \begin{subfigure}{0.32\textwidth}
    \includegraphics[width=\textwidth]{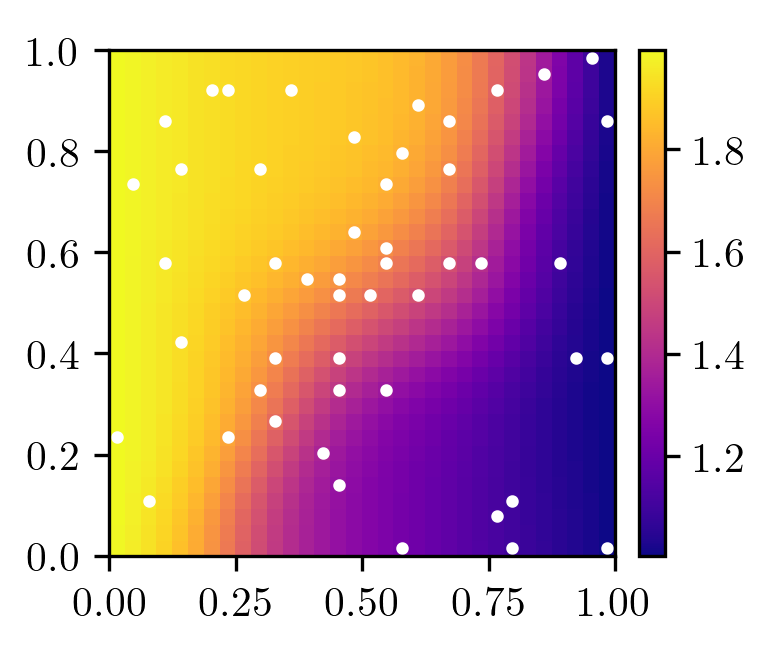}
    \caption{$u(x)$, Gaussian kernel}
  \end{subfigure}
  \begin{subfigure}{0.32\textwidth}
    \includegraphics[width=\textwidth]{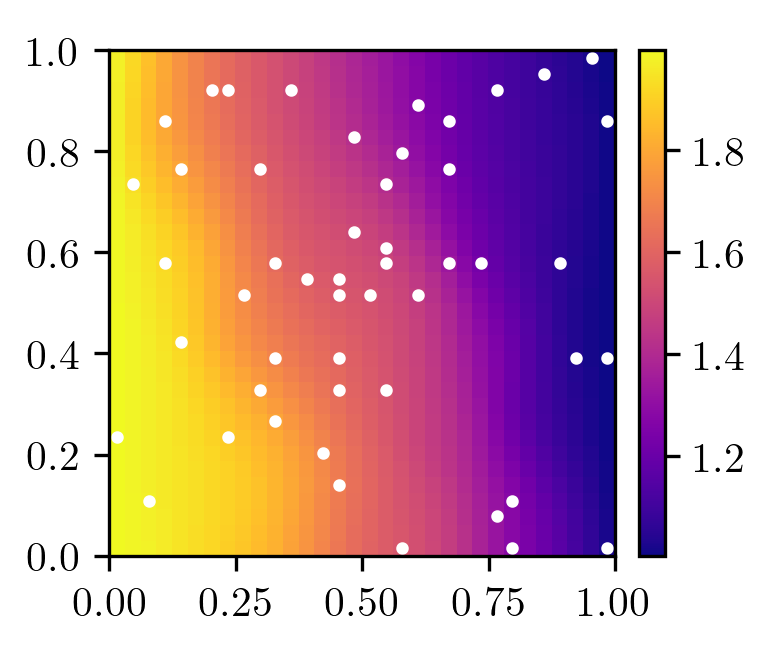}
    \caption{$u(x)$, Mat\'{e}rn $\nu = 5 / 2$}
  \end{subfigure}
  \begin{subfigure}{0.32\textwidth}
    \includegraphics[width=\textwidth]{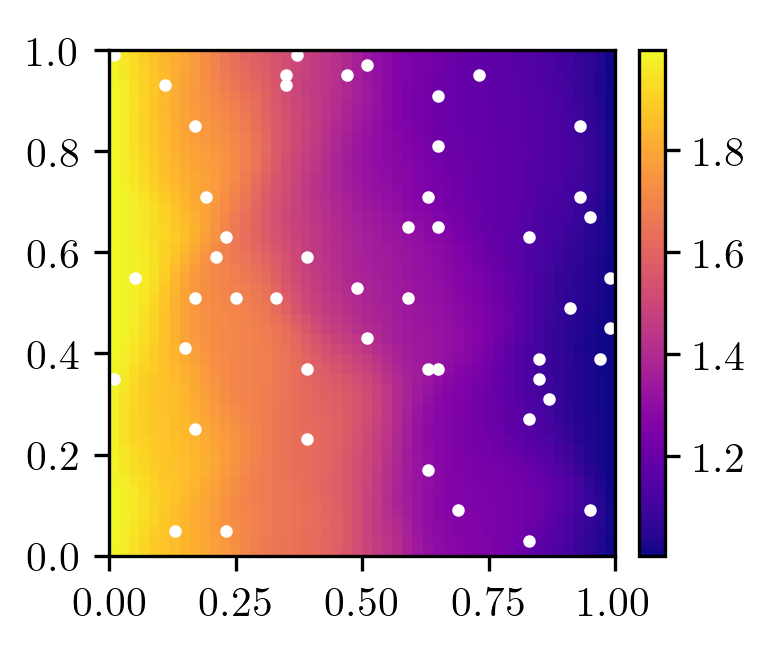}
    \caption{$u(x)$, Mat\'{e}rn $\nu = 3 / 2$}
  \end{subfigure}
  \caption{Sample reference $y$ (above) and $u$ (below) fields for the parameters of \cref{tbl:numerical-continuous-params}.}
  \label{fig:Yref-uref-02-50}
\end{figure}

We first consider continuous reference $y$ fields of various degrees of smoothness.
Three reference fields are generated as realizations of zero-mean Gaussian processes with the isotropic Mat\'{e}rn (with $\nu = \{ 3/ 2, 5/2\}$) and Gaussian kernels for the values of the kernel hyperparameters (namely the correlation length $\lambda$ and standard deviation $\sigma$), listed in~\cref{tbl:numerical-continuous-params}.
The corresponding reference $u$ fields are computed by solving the FV discretization of the PDE of \cref{eq:numerical-pde}.
Finally, observation locations for the $y$ and $u$ are chosen randomly from the set of FV cell centers.
The number of observations are listed in~\cref{tbl:numerical-continuous-params}.
For the reference fields with the Gaussian and  Mat\'{e}rn $\nu = 5/2$ kernels, we assume that $50$ observations for both  $y$ and $u$ are available. For the Mat\'{e}rn $\nu = 3/2$ case, we use $200$ $y$ observations to estimate the field.
The larger number of observations is necessary for the latter case as this problem is more challenging due to its short correlation length and lower smoothness.
The number of KLE terms for this kernel given by the condition~\cref{eq:truncation} with $\text{rtol} = 99\%$ is $511$.
The reference fields and observation locations are shown in~\cref{fig:Yref-uref-02-50}.

As described in~\cref{sec:ckli-ckle-y}, we construct the cKLE model for $y$ by training a GPR model to the observation data $(\mathbf{y}_{\mathrm{s}}, X^y_{\mathrm{s}})$.
We consider two scenarios: (i) the data-generating kernel and its hyperparameters are known, and (ii) the data-generating kernel is known but the hyperparameters are unknown.
We will refer to the first scenario as ``cKLI'' and the second scenario as ``cKLI-$\theta$''.
For cKLI-$\theta$, the kernel hyperparameters are estimated using the observations $y_{\mathrm{s}}$ via marginal likelihood estimation, which is performed using the library \textsf{GPy}~\citep{gpy_2014}.

\begin{figure}[tbhp]
  \centering
  \begin{subfigure}{0.32\textwidth}
    \includegraphics[width=\textwidth]{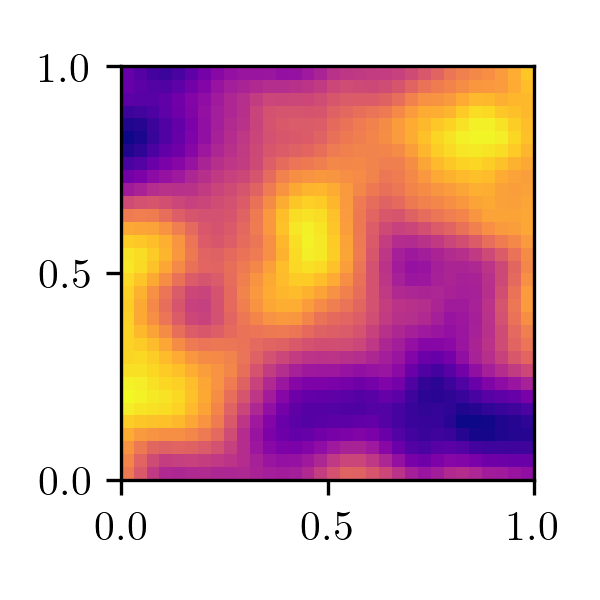}
    \caption{Mat\'{e}rn $\nu = 5/2$ ref.}
  \end{subfigure}
  \begin{subfigure}{0.32\textwidth}
    \includegraphics[width=\textwidth]{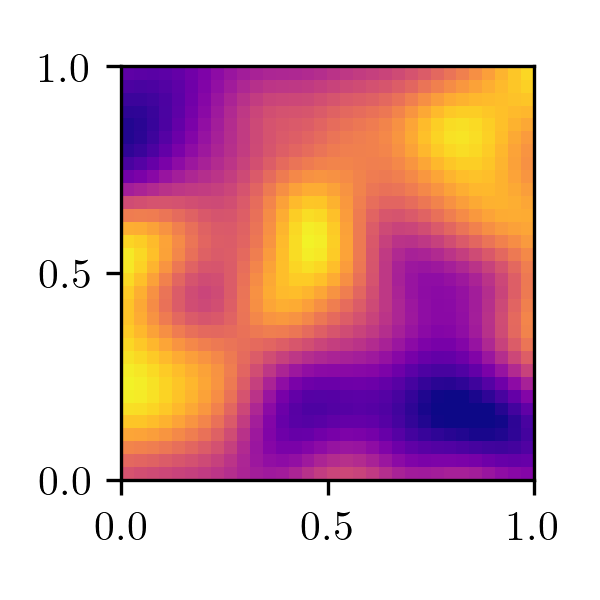}
    \caption{cKLI-$\theta$}
  \end{subfigure}
  \begin{subfigure}{0.32\textwidth}
    \includegraphics[width=\textwidth]{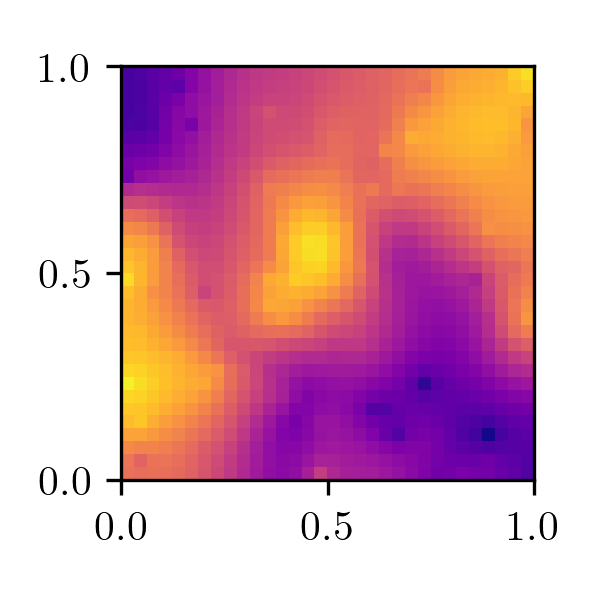}
    \caption{MAP}
  \end{subfigure}\\
  \begin{subfigure}{0.32\textwidth}
    \includegraphics[width=\textwidth]{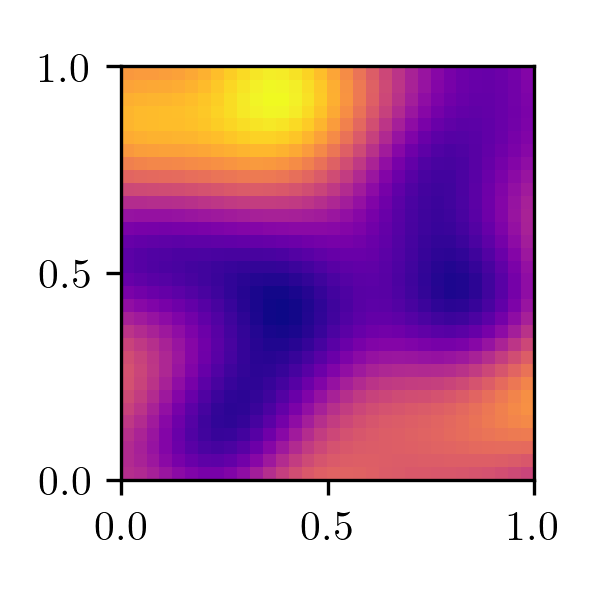}
    \caption{Gaussian kernel ref.}
  \end{subfigure}
  \begin{subfigure}{0.32\textwidth}
    \includegraphics[width=\textwidth]{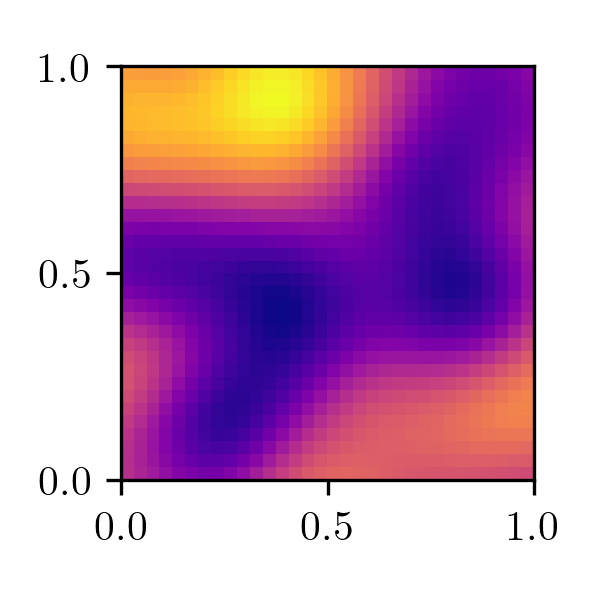}
    \caption{cKLI-$\theta$}
  \end{subfigure}
  \begin{subfigure}{0.32\textwidth}
    \includegraphics[width=\textwidth]{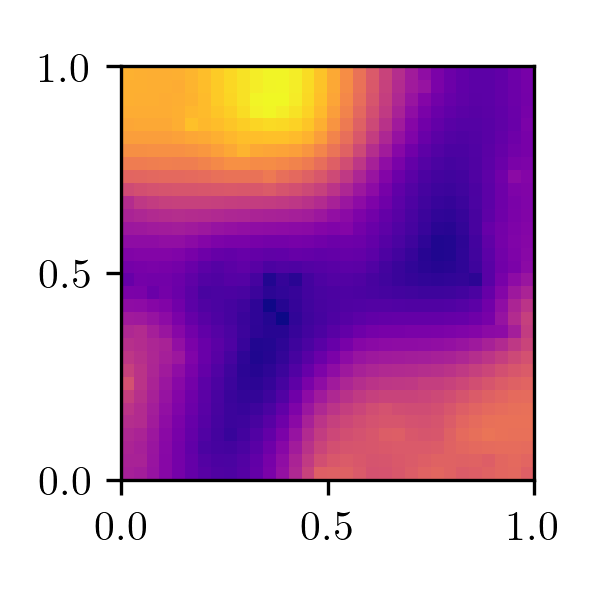}
    \caption{MAP}
  \end{subfigure}\\
  \begin{subfigure}{0.32\textwidth}
    \includegraphics[width=\textwidth]{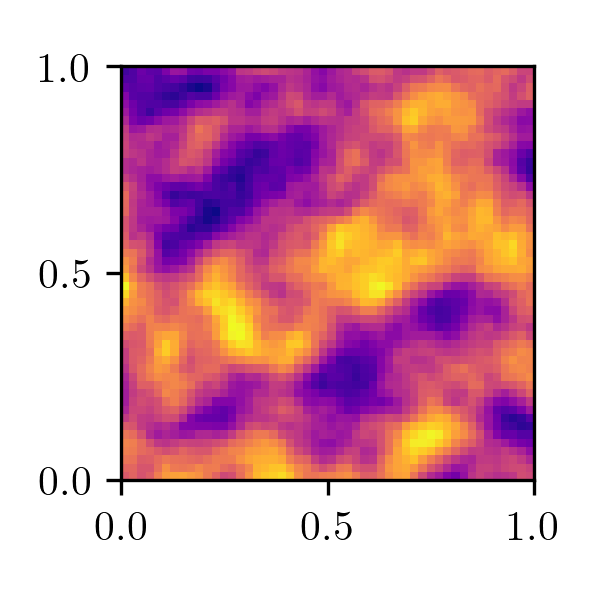}
    \caption{Mat\'{e}rn $\nu = 3/2$ ref.}
  \end{subfigure}
  \begin{subfigure}{0.32\textwidth}
    \includegraphics[width=\textwidth]{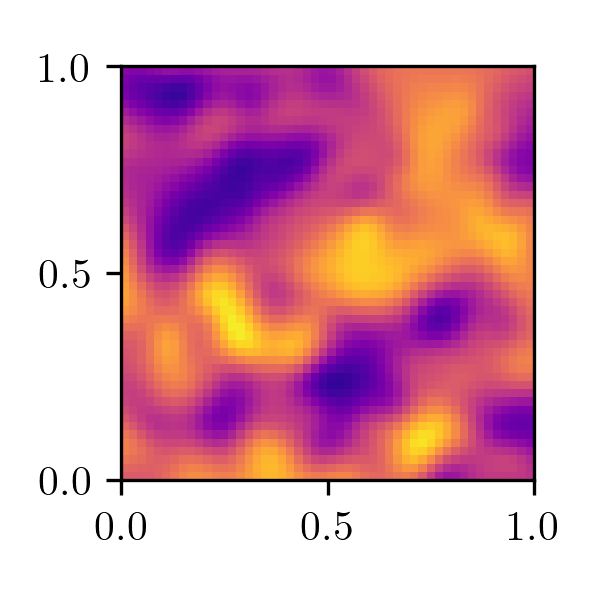}
    \caption{cKLI-$\theta$}
  \end{subfigure}
  \begin{subfigure}{0.32\textwidth}
    \includegraphics[width=\textwidth]{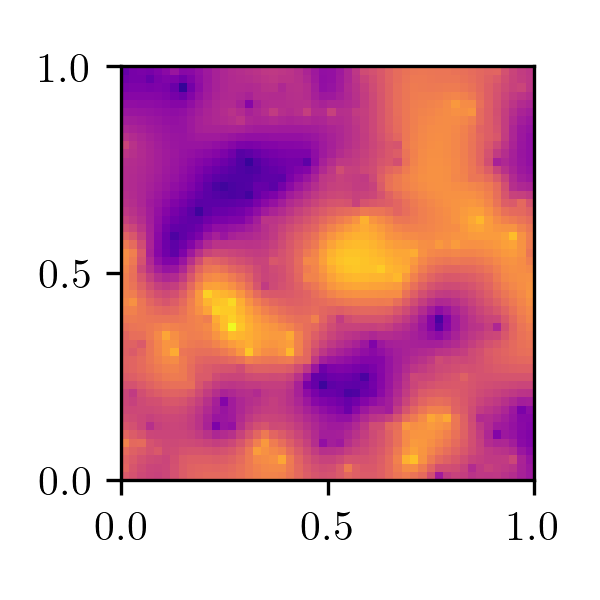}
    \caption{MAP}
  \end{subfigure}
  \caption{cKLI-$\theta$ estimate (middle) and MAP estimate (right) of the reference log-diffusion fields of \cref{fig:Yref-uref-02-50}.}
  \label{fig:Yest-02-50}
\end{figure}

\cref{fig:Yest-02-50} presents the reference $y$ fields and the cKLI-$\theta$ and MAP estimates of these fields.
The cKLI-$\theta$ estimates are computed using $N_{\xi}$ and $N_{\eta}$ listed in \cref{tbl:numerical-continuous-params}.
For all cases, we used $N_{\mathrm{ens}} = \num{5e3}$.
It can be seen that the cKLI-$\theta$ estimates of $y$ are more accurate than the MAP estimates for all considered cases; this advantage is more noticeable for the Mat\'{e}rn cases, which is less smooth than the Gaussian case.
Furthermore, the MAP estimates exhibit peaks at the $y$ observation locations, a phenomenon typical to $H_1$ regularization, whereas the cKLI-$\theta$ estimates are smooth and thus better approximate the reference fields.

We define the ``relative $\ell_p$ error'' as the $\ell_p$-norm of the estimation error with respect to the $\ell_p$-norm of the reference field, that is,
\begin{equation*}
  \text{relative $\ell_p$ error} \coloneqq \| y_{\mathrm{ref}} - y_{\mathrm{est}} \|_p  / \| y_{\mathrm{ref}} \|_p.
\end{equation*}%
In \cref{tbl:rel-l2-error-continuous}, we present the relative $\ell_2$ error of the estimates shown in~\cref{fig:Yest-02-50}.
We also present the relative error of the cKLI-$\theta$ estimate computed using subsampled residuals with a subsampling factor of $2$ in each direction (resulting in a reduction in the dimension of the vector of residuals by a factor of $4$).
Subsampling reduces the dimension of the vector of residuals and therefore reduces the computational effort of computing the cKLI-$\theta$ estimate.
For comparison, we also present the $\ell_2$ error in the estimates obtained with the cKLI method.
 As expected, the accuracy of the cKLI estimated $y$ is the same or better than that of the cKLI-$\theta$ estimate for all considered cases.
This is because the accuracy of PICKLE estimation depends on the accuracy of the estimated $y$ kernel, and in the cKLI case we assume that the $y$ kernel is known exactly.
It can be seen that for the cases considered so far, the PICKLE for estimating $y$ is more accurate in the $\ell_2$ sense than MAP, and that subsampling by a factor of $2$ of the vector of residuals does not significantly increase the PICKLE estimation error.

For comparison, we also estimate $y$ using the PINN method for parameter estimation~\cite{tartakovsky2018learning}.
Here, we represent the $y$ and $u$ fields using deep feed-forward neural networks with three hidden layers and 30 neurons per layer.
The residual is estimated at $N_r = 1024$ points.
We conduct ten simulations with different initializations using the Xavier's initialization scheme, and train the PINN networks by using the L-BFGS-B method.
The mean and standard deviation across initializations of the relative $\ell_2$ error is reported in \cref{tbl:rel-l2-error-continuous}.
The relative $\ell_2$ error of PICKLE cKLI-$\theta$ estimation is approximately 320\% smaller than of PINN estimation for the Gaussian kernel reference, 70\% smaller for the Mat\'{e}rn $\nu = 5/2$ reference, and 50\% for the Mat\'{e}rn $\nu = 5/2$ reference.

\begin{table}
  \centering
    \caption{Relative $\ell_2$ error of the $y$ estimates shown in \cref{fig:Yest-02-50} and obtained with PINNs.
  For cKLI-$\theta$, ``Full'' indicates the estimate computed using the full vector of FV residuals, and ``Subsampled'' indicates the estimate computed using a subsampling of the vector of residuals by a factor of $2$ in each spatial direction.}
  \begin{tabular}{lccccccc}
    & \multicolumn{2}{c}{cKLI} & \multicolumn{2}{c}{cKLI-$\theta$} & MAP & PINN \\
              & Full & Subsampled & Full & Subsampled & & \\
    \midrule
    Gaussian & 0.010 & 0.015 & 0.017 & 0.022 & 0.150 & 0.072(27)\\
    Mat\'{e}rn $\nu = 5/2$ & 0.099 & 0.109 & 0.099 & 0.111 & 0.224 & 0.169(11)\\
    Mat\'{e}rn $\nu = 3/2$ & 0.257 & 0.254 & 0.262 & 0.261 & 0.419 & 0.388(16)
  \end{tabular}
  \label{tbl:rel-l2-error-continuous}
\end{table}

To evaluate the robustness of PICKLE, we calculate the relative $\ell_2$ estimation error for different reference fields (generated as realizations of the random fields with the Mat\'{e}rn  ($\nu = 5/2$) and Gaussian kernels with $\sigma = 1.0$ and the correlation lengths $\lambda = 0.2$ and $0.5$.
) and choices of observation locations.
Furthermore, we study how the relative $\ell_2$ error depends on the number of cKLE terms in the expansions of the $y$ and $u$ fields.
For each combination of kernel and correlation length, we generate 10 reference $y$ fields and the corresponding $u$ fields.
For each reference field, we randomly generate observation locations, and compute the cKLI and cKLI-$\theta$ estimates of $y$.
We do this for $N^u_{\mathrm{s}} = 50$, $N_{\eta} = 100$, two values of $N^y_{\mathrm{s}}$, $10$ and $50$, and various values of $N_{\xi}$.

\begin{figure}[tbhp]
  \centering
  \begin{subfigure}{0.48\textwidth}
    \includegraphics[width=\textwidth]{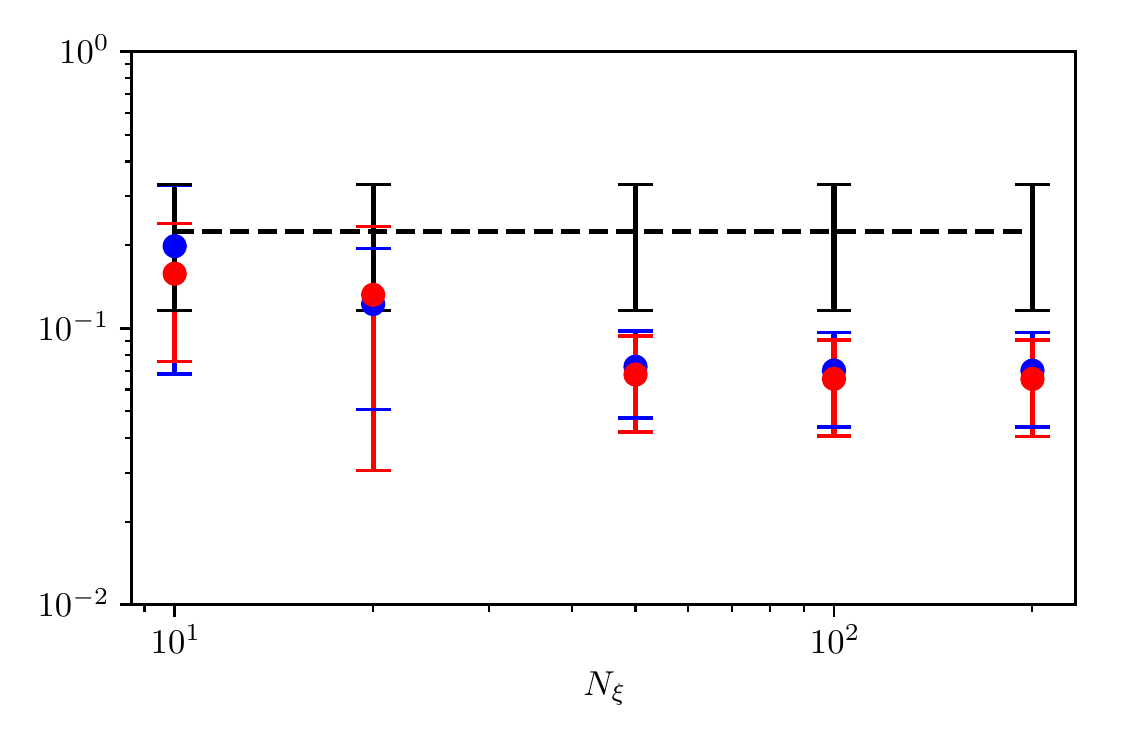}
    \caption{$\lambda = 0.5$, $N^y_{\mathrm{s}} = 10$}
  \end{subfigure}
  \begin{subfigure}{0.48\textwidth}
    \includegraphics[width=\textwidth]{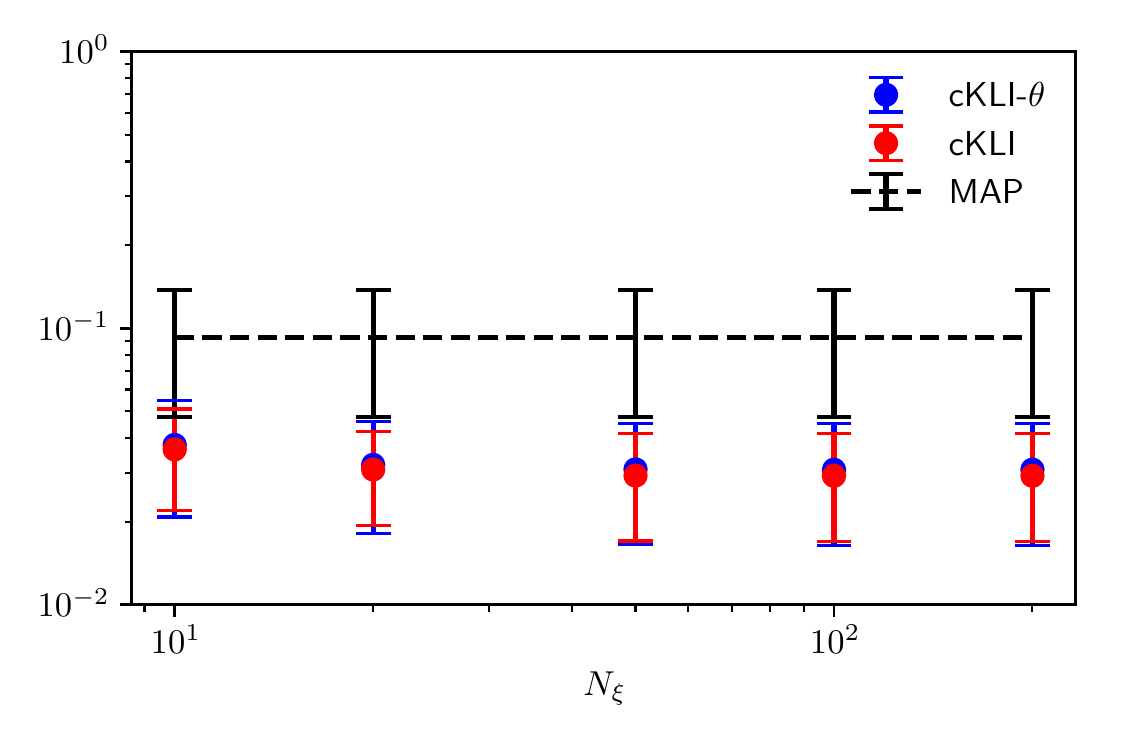}
    \caption{$\lambda = 0.5$, $N^y_{\mathrm{s}} = 50$}
  \end{subfigure}\\
  \begin{subfigure}{0.48\textwidth}
    \includegraphics[width=\textwidth]{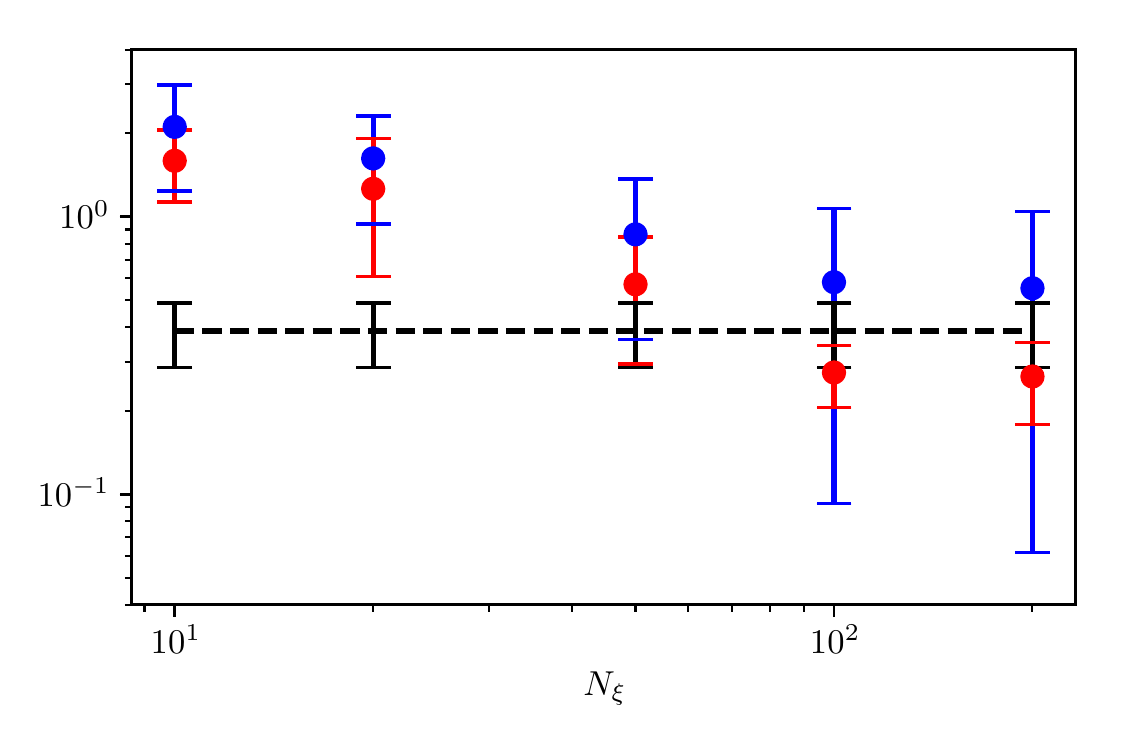}
    \caption{$\lambda = 0.2$, $N^y_{\mathrm{s}} = 10$}
    \label{fig:convergence-wrt-NYxi-matern-bad}
  \end{subfigure}
  \begin{subfigure}{0.48\textwidth}
    \includegraphics[width=\textwidth]{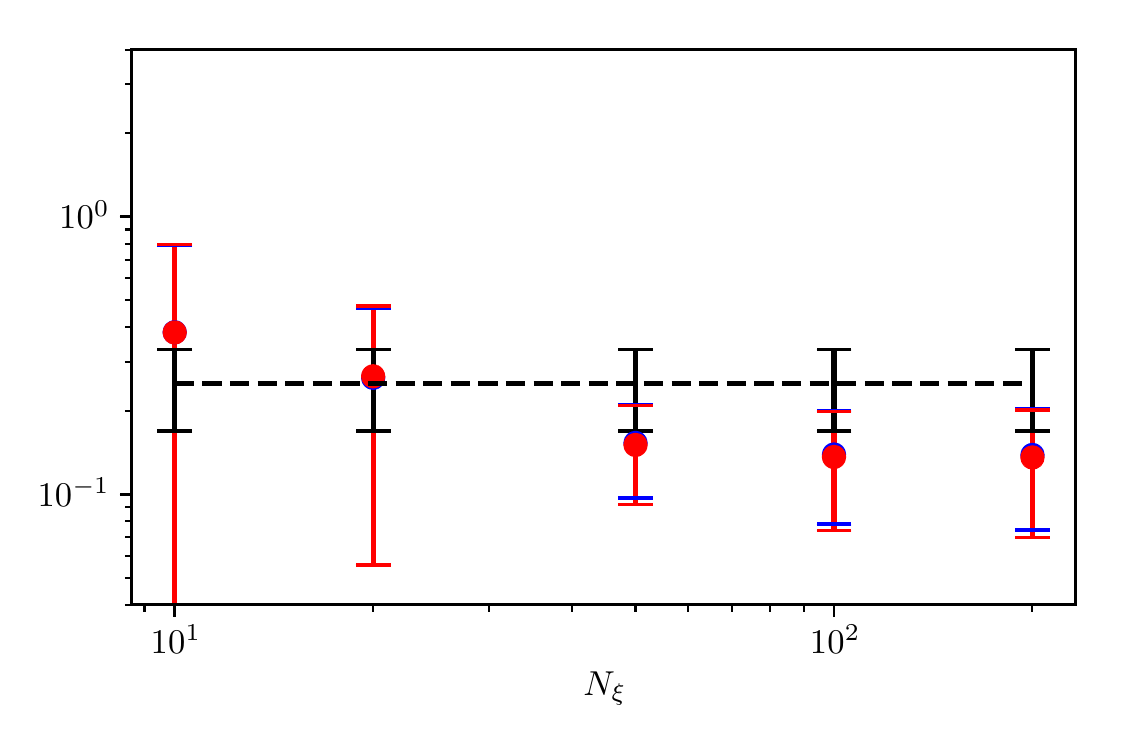}
    \caption{$\lambda = 0.2$, $N^y_{\mathrm{s}} = 50$}
    \label{fig:convergence-wrt-NYxi-matern-good}
  \end{subfigure}
  \caption{Relative $\ell_2$ error for Mat\'{e}rn covariance kernel with $\nu = 5/2$ and with  different values of $\lambda$ and $N^y_{\mathrm{s}}$}
  \label{fig:convergence-wrt-NYxi-matern}
\end{figure}

\begin{figure}[tbhp]
  \centering
  \begin{subfigure}{0.48\textwidth}
    \includegraphics[width=\textwidth]{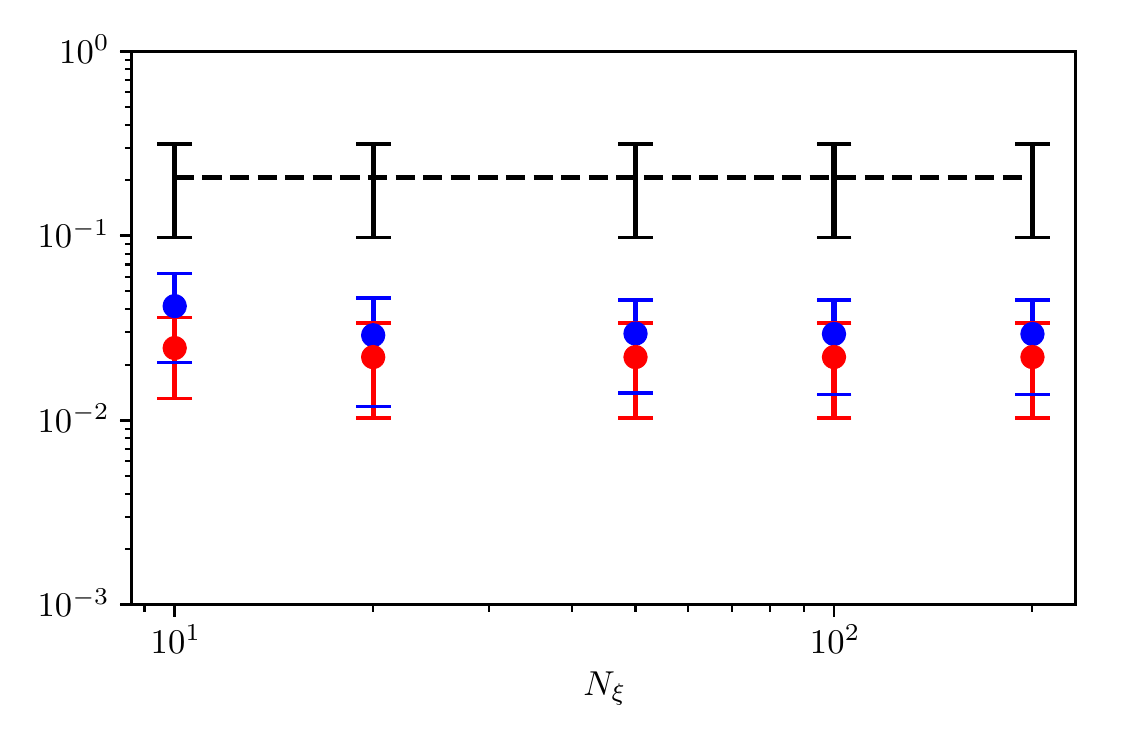}
    \caption{$\lambda = 0.5$, $N^y_{\mathrm{s}} = 10$}
  \end{subfigure}
  \begin{subfigure}{0.48\textwidth}
    \includegraphics[width=\textwidth]{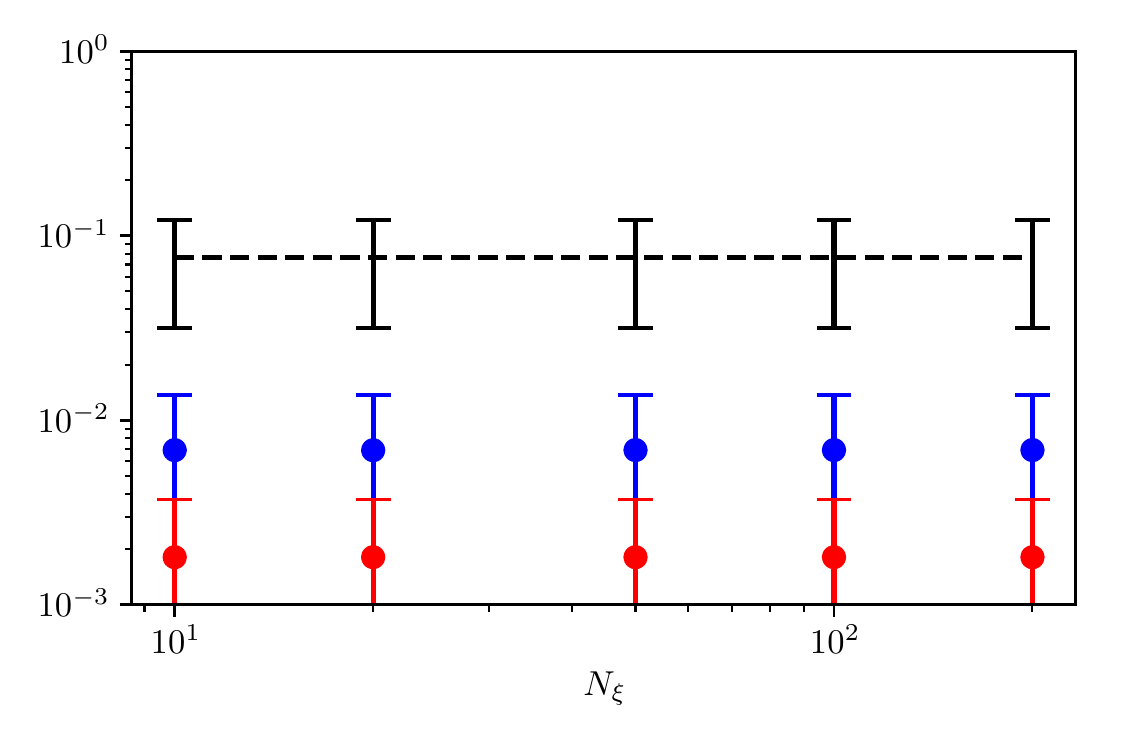}
    \caption{$\lambda = 0.5$, $N^y_{\mathrm{s}} = 50$}
  \end{subfigure}\\
  \begin{subfigure}{0.48\textwidth}
    \includegraphics[width=\textwidth]{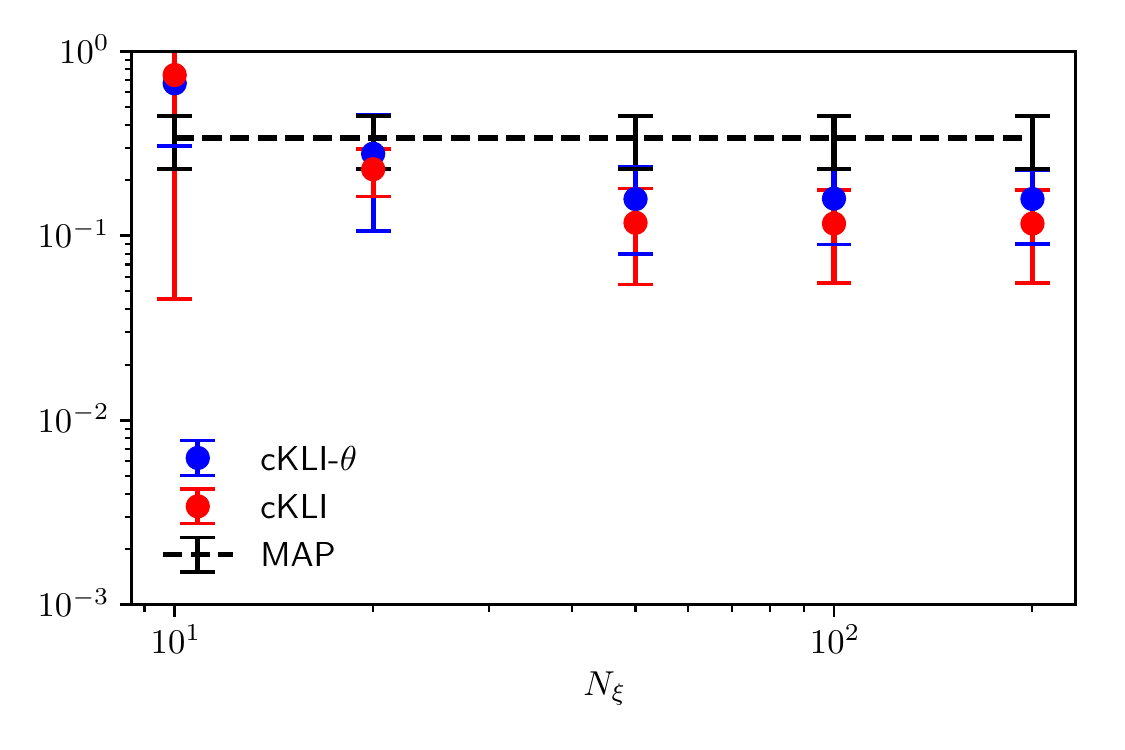}
    \caption{$\lambda = 0.2$, $N^y_{\mathrm{s}} = 10$}
  \end{subfigure}
  \begin{subfigure}{0.48\textwidth}
    \includegraphics[width=\textwidth]{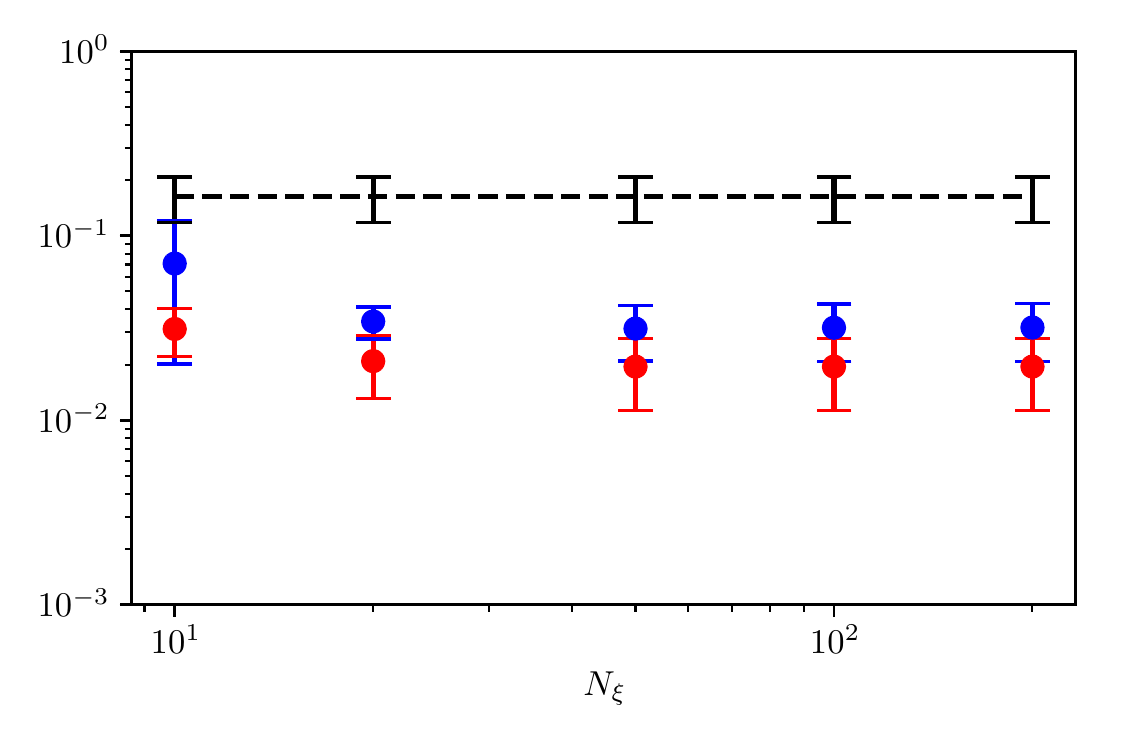}
    \caption{$\lambda = 0.2$, $N^y_{\mathrm{s}} = 50$}
  \end{subfigure}
  \caption{Relative $\ell_2$ error for the Gaussian covariance kernel with different values of $\lambda$ and $N^y_{\mathrm{s}}$.}
  \label{fig:convergence-wrt-NYxi-gaussian}
\end{figure}

The relative $\ell_2$ error as a function of $N_{\xi}$ is shown in \cref{fig:convergence-wrt-NYxi-matern,fig:convergence-wrt-NYxi-gaussian}.
As in the previous example, we can see that the accuracy of the cKLI estimated $y$ is the same or better than that of the cKLI-$\theta$ estimate for all considered cases.
The cKLI estimates are consistently more accurate in the $\ell_2$ sense than the MAP estimate for sufficiently large $N_{\xi}$.
In particular, we note that a good rule of thumb for $N_{\xi}$ is to be larger than the number of KLE terms of the reference kernel for $\text{rtol = 99\%}$ minus the number of observations.

As expected, for the rougher Mat\'{e}rn kernel, more KL terms are needed to obtain an accurate $y$ estimate than for the smoother Gaussian kernel.
The same observation is true with respect to the correlation length: The smaller is the correlation length, the more KL terms are needed to obtain an accurate $y$ estimate. 

The cKLI-$\theta$ estimates require additional discussion.
For all considered fields except one, the cKLI-$\theta$ estimate of $y$ is more accurate than the MAP estimate for sufficiently large $N_\xi$.
For $N^y_{\mathrm{s}} = 10$ and the rough (Mat\'{e}rn) kernel with small correlation length, the cKLI-$\theta$ is worse than the MAP estimate for all considered  $N_{\xi}$ (\cref{fig:convergence-wrt-NYxi-matern-bad}). 
This is because $10$ $y$ observations are not sufficient to obtain adequate estimates of the hyperparameters of the kernel for such a rough $y$ field.
\cref{fig:convergence-wrt-NYxi-matern-good} shows that for the same Mat\'{e}rn kernel, a very accurate estimate of hypermaparameters is obtained with 50 $y$ measurements, and the cKLI-$\theta$ estimates of $y$ are as accurate as the cKLI estimates and are more accurate than MAP estimation for sufficiently large values of $N_\xi$.  

The comparison of the cKLI and cKLI-$\theta$ results show that the more accurate estimate of the $y$ kernel is available, the less terms in the cKLE model of $y$ are needed to obtain an accurate estimate of $y$.
These results also indicate that it is not necessary to know the $y$ kernel exactly for PICKLE estimation to produce an accurate estimate of $y$, given that $N_{\xi}$ is sufficiently large.    

\subsection{Piecewise-constant diffusion field}
\label{sec:numerical-piecewise}

\begin{figure}[tbhp]
  
  \sbox\threesubbox{%
    \resizebox{\dimexpr.96\textwidth}{!}{%
      \includegraphics[height=3cm]{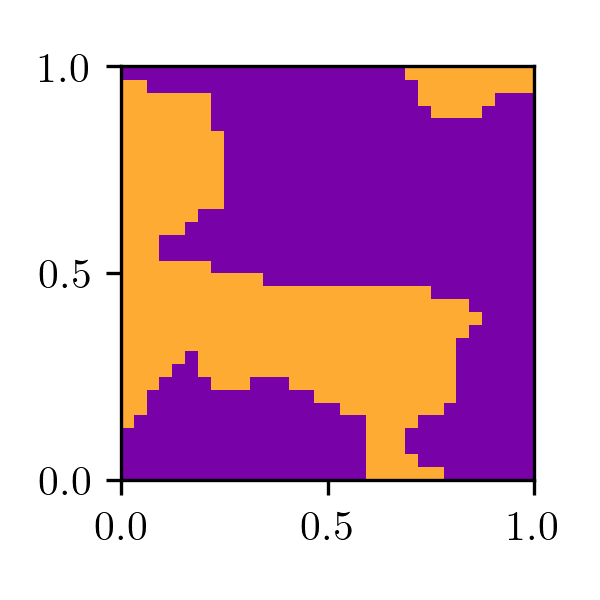}%
      \includegraphics[height=3cm]{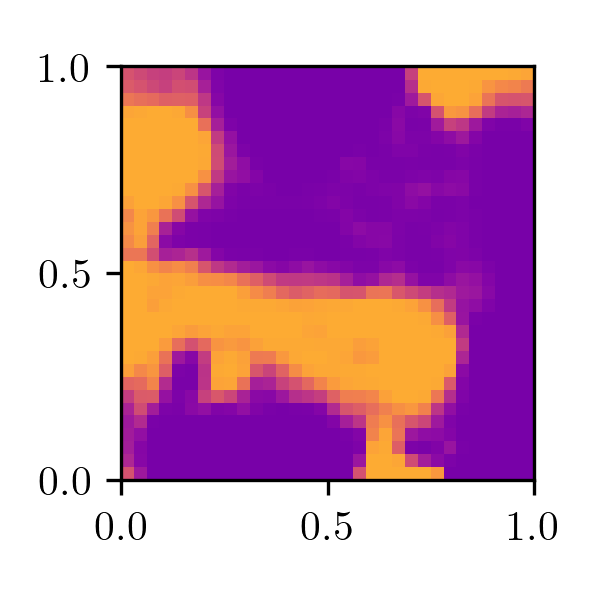}%
      \includegraphics[height=3cm]{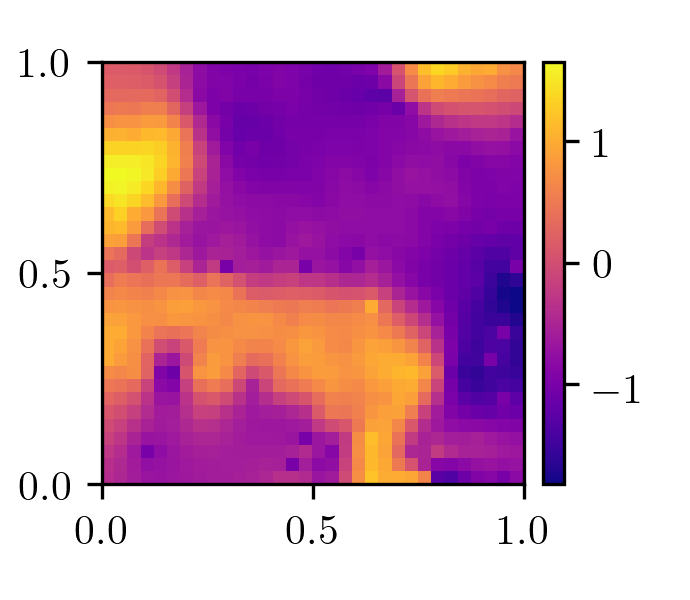}
    }%
  }%
  \setlength{\threesubht}{\ht\threesubbox}

  \centering
  \begin{subfigure}{0.32\textwidth}
    \includegraphics[height=\threesubht]{figures/Yref_binary}
    \caption{Ref.}
  \end{subfigure}
  \begin{subfigure}{0.32\textwidth}
    \includegraphics[height=\threesubht]{figures/YCKLI_binary}
    \caption{cKLI-$\theta$}
  \end{subfigure}
  \begin{subfigure}{0.32\textwidth}
    \includegraphics[height=\threesubht]{figures/YMAP_binary}
    \caption{MAP}
  \end{subfigure}
  
  \caption{Reference piecewise-continuous log-diffusion field, and estimates computed usign cKLI-$\theta$ and MAP.
    Relative $\ell_1$ error of cKLI-$\theta$ is $0.179$ and of MAP is $0.380$.}
  \label{fig:Yest-binary}
\end{figure}

Finally, we consider the case of a piecewise-constant diffusion field $k$ (\cref{fig:Yest-binary}).
Due to Gibb's phenomenon, a large number of KLE terms would be necessary to accurately represent either $k$ or $y \coloneqq \log k$ directly~\cite{tagade2017mitigating}.
Therefore, to treat this case with our proposed method, we introduce a latent field $f(x)$ that can be represented accurately using a finite-dimensional representation.
For this application we assume that the log-diffusion field consist of two facies, with constant log-diffusion values of $y_1$ and $y_2$, with $y_1 > y_2$.
The log-diffusion coefficient is then approximated in terms of the latent field $f(x)$ as
\begin{equation}
  \label{eq:log-k-binary}
  y(x) = (y_1 - y_2) \operatorname{expit} \left ( \varepsilon^{-1} f(x) \right ) + y_2,
\end{equation}
where $\operatorname{expit} \coloneqq 1 / [1 + \exp (-x)]$ is the logistic function, and $\varepsilon > 0$ is a small constant~\cite{menard2002applied}.
In the limit $\varepsilon \to 0$, $\operatorname{expit} \left ( \varepsilon^{-1} ( \cdot) \right )$ approximates the step function from $0$ to $1$.

To compute the PICKE cKLI-$\theta$ estimate of $y$, we first construct a cKLE for the latent field $f$ from sparse measurements of $y$, which we accomplish via GP classification~\citep{williams2006gaussian}.
As the latent field is not observed directly, we cannot use the GPR~\cref{eq:gpr-mean,eq:gpr-var} to construct the conditional GP model $\hat{f}^c(x, \omega)$.
Instead, we proceed as follows: The observations $\mathbf{y}_{\mathrm{s}} = (y_1, \dots, y_{N^{y}_{\mathrm{s}}})^{\top}$ are translated into a vector of binary values $\mathbf{b}_{\mathrm{s}} = (b_i, \dots, b_{N^y_{\mathrm{s}}})^{\top}$, where $b_i = 0$ and $b_i = 1$ indicate $y_i=y_2$ and $y_i=y_1$, respectively.
These binary observations are employed to construct the logistic GP classifier $\hat{f}^c(x, \omega)$, corresponding to the random field $\hat{f}$ conditioned on the outcomes $b_i$ of the Bernoulli random variables (e.g. random variables with binary outcome) with probability of $b = 1$ given by $\operatorname{expit}(\hat{f}(X^y_i))$\footnote{%
  Note that GPR, \cref{eq:gpr-mean,eq:gpr-var}, can be understood in the same terms.
  Specifically, $\hat{y}^c$ is equivalent to the random field $\hat{f}$ conditioned on the outcomes $\mathbf{y}_{\mathrm{s}}$ of the random variable $\mathcal{N}(\hat{f}(X^y_{\mathrm{s}}), \Sigma)$.
}, that is,
\begin{equation*}
  b_i \sim \operatorname{Bernoulli} ( \operatorname{expit} ( \hat{f}(X^y_i) ) ), \quad i \in [1, N^y_{\mathrm{s}}].
\end{equation*}
The conditioning is performed using the expectation propagation algorithm as implemented by the library \textsf{GPy}~\citep{gpy_2014}.
Once the conditional mean and covariance have been estimated, we then compute the cKLE of $\hat{f}^c$.

Next, we construct the sampling-based covariance model for $u$ by using~\cref{alg:uc}.
The realizations $\{y^{(i)} \}$ are generated by sampling fields $\{ f^{(i)} \}$ from the cKLE of $\hat{f}^c$, which are then substituted into~\cref{eq:log-k-binary}.
Once the conditional covariance of $u$ is found, the PICKLE estimate of $f$ (and of $y$ through~\cref{eq:log-k-binary}) is computed using \cref{alg:ckli}.

\cref{fig:Yest-binary} shows the reference binary $y$ field and the PICKLE cKLI-$\theta$ and MAP estimates of $y$ using $25$ measurements of $y$ and $100$ measurements of $u$.
The reference $f$ field is generated as a realization of the zero-mean Gaussian process with the isotropic Mat\'{e}rn ($\nu = 5/2$) kernel, $\sigma = 1.0$, and $\lambda = 0.2$.
The reference $y$ field is generated by substituting the reference $f$ field into~\cref{eq:log-k-binary} with $\epsilon = 100$.
As before, the reference $u$ field is generated by solving~\cref{eq:numerical-pde} for the reference $y$ field.
It can be seen that the PICKLE cKLI-$\theta$ estimate of $y$ is closer to the reference $y$ and has a significantly sharper boundary between the ``$y_1$'' and ``$y_2$'' regions than the MAP estimated $y$.
The relative $l_1$ error of the PICKE cKLI-$\theta$ estimate of $y$ is 0.179, more than two times smaller than the MAP estimation error of 0.380.
These results indicate that PICKLE estimation can be employed to estimate discontinuous fields by expressing these fields in terms of cKLEs of continuous latent fields.

\section{Conclusions}
\label{sec:conclusions}

We presented a new physics-informed machine learning approach, termed PICKLE, for learning parameters and states of stationary physical systems from sparse measurements constrained by the stationary PDE models governing the behavior of said systems.
In PICKLE, parameters and states are approximated using cKLEs, i.e., KLEs conditioned on measurements, resulting in low-dimensional models of spatial fields that honor observed data. 
Finally, the coefficients in the cKLEs are estimated by minimizing the norm of the residual of the PDE model evaluated at a finite set of points in the computational domain, ensuring that the reconstructed parameters and states are consistent with both the observations and the PDE model to an arbitrary level of accuracy.

The cKLEs are constructed using the eigendecomposition of covariance models of spatial variability. 
For the model parameter (space-dependent diffusion coefficient), we employed a parameterized covariance model calibrated on parameter observations; for the model state, the covariance was estimated from a number of forward simulations of the PDE model corresponding to realizations of the parameter drawn from its  cKLE.
We demonstrated that the accuracy of the PICKLE method depends on the accuracy of the estimated parameter covariance, which in turn depends on the number of measurements.
It is important to note that transfer learning could be used to estimate the covariance of parameters, e.g., measurements collected in other systems with statistically similar properties can be used to estimate the covariance function of the model parameters.

We applied PICKLE to solve an inverse problem associated with the steady-state diffusion equation with unknown space-dependent diffusion coefficient. Specifically, we  used PICKLE to estimate the log-diffusion coefficient from sparse measurements of the log-diffusion coefficient and the state of the system.
We considered continuous and discontinuous diffusion coefficients. 
For continuous diffusion coefficients with different degrees of roughness (corresponding to different covariance kernels and correlation lengths), we demonstrated that the PICKLE estimates of the  diffusion coefficient are more accurate than those of the MAP and physics-informed neural networks (PINN) method.
The comparison with the PINN method suggests that cKLEs are better representations of sparsely-measured spatially-correlated fields than neural networks. We also found that PICKLE provides a better estimate of the discontinuous conductivity field than the MAP method.

Our results indicate that the PICKLE method can be used for estimating space-dependent parameters and states regardless of their underlying statistical distribution. Even though the cKLE expansion in PICKLE is constructed using the GPR estimates of the mean and covariance functions, we demonstrated that accurate estimates can be obtained when cKLE is used to model fields with highly non-Gaussian statistics, including the solution of the diffusion equation on the bounded domain and the discontinuous diffusion coefficient.   

  

\section*{Acknowledgments}

This work was supported by the Applied Mathematics Program within the U.S. Department of Energy Office of Advanced Scientific Computing Research.
Pacific Northwest National Laboratory is operated by Battelle for the DOE under Contract DE-AC05-76RL01830.

\bibliographystyle{elsarticle-num}
\bibliography{conditioned_kl_inversion}

\end{document}